\documentclass[10pt,letterpaper]{amsart}
\usepackage{charter}
\usepackage[OT2,OT1]{fontenc}
\usepackage[latin9]{inputenc}
\usepackage{array}
\usepackage{calc}
\usepackage{units}
\usepackage{textcomp}
\usepackage{mathrsfs}
\usepackage{bm}
\usepackage{multirow}
\usepackage{amstext}
\usepackage{amsthm}
\usepackage{amssymb}
\usepackage{cancel}
\usepackage{stmaryrd}
\usepackage{graphicx}
\usepackage{rotfloat}

\makeatletter

\pdfpageheight\paperheight
\pdfpagewidth\paperwidth

\numberwithin{equation}{section}
\numberwithin{figure}{section}
\theoremstyle{definition}
    \ifx\thechapter\undefined
      \newtheorem{defn}{\protect\definitionname}
    \else
      \newtheorem{defn}{\protect\definitionname}[chapter]
    \fi
\theoremstyle{plain}
    \ifx\thechapter\undefined
	    \newtheorem{thm}{\protect\theoremname}
	  \else
      \newtheorem{thm}{\protect\theoremname}[chapter]
    \fi
\theoremstyle{remark}
    \ifx\thechapter\undefined
      \newtheorem{claim}{\protect\claimname}
    \else
      \newtheorem{claim}{\protect\claimname}[chapter]
    \fi
\theoremstyle{definition}
    \ifx\thechapter\undefined
      \newtheorem{example}{\protect\examplename}
    \else
      \newtheorem{example}{\protect\examplename}[chapter]
    \fi
\theoremstyle{plain}
    \ifx\thechapter\undefined
      \newtheorem{prop}{\protect\propositionname}
    \else
      \newtheorem{prop}{\protect\propositionname}[chapter]
    \fi

\usepackage{amsfonts}
\usepackage{mathrsfs}
\usepackage{graphicx}
\usepackage{amscd}
\usepackage{enumerate}
\usepackage{url}

\newenvironment{acknowledgements}{%
  \section*{Acknowledgements}%
}{}
\newcommand{\refeq}[1]{\eqref{#1}}
 \newcommand{\cyr}{
\renewcommand\rmdefault{wncyr} \renewcommand\sfdefault{wncyss} \renewcommand\encodingdefault{OT2} \normalfont
\selectfont
}
\DeclareTextFontCommand{\textcyr}{\cyr}    

   \def\@settitle
{
\begin{center} \baselineskip14\p@\relax \LARGE\@title \end{center}
}

\usepackage[normalem]{ulem}
\usepackage[mathscr]{eucal}
\numberwithin{equation}{section}

\usepackage[table]{xcolor}
\usepackage{verbatim}
\usepackage{amscd}
\usepackage{bm}\usepackage[mathscr]{eucal}
\usepackage[all]{xy}
\usepackage{enumerate}
\usepackage{color}
\usepackage{colortbl}\usepackage{comment}
\usepackage{here}
\usepackage{rotating} 
\usepackage[table]{xcolor} 
\usepackage{pdflscape}

\input xyimport.tex

\newenvironment{spmatrix}{\left ( \begin{smallmatrix}}{\end{smallmatrix}\right)}

\title{\bf Constructing prime $\mQ$-Fano threefolds \\ of codimension four via key varieties\\ related with $\mP^2\times \mP^2$-fibrations}
\author{Hiromichi Takagi}

\address{Department of Mathematics, Gakushuin University, 1-5-1, Mejiro, Toshima, Tokyo, 171-8588 Japan}
\email{hiromici@math.gakushuin.ac.jp}

\newcommand{\sO}{\mathcal{O}}

\newcommand{\mA}{\mathbb{A}}
\newcommand{\mC}{\mathbb{C}}

\newcommand{\mN}{\mathbb{N}}
\newcommand{\mP}{\mathbb{P}}
\newcommand{\mQ}{\mathbb{Q}}
\newcommand{\mZ}{\mathbb{Z}}

\newcommand{\Bs}{\mathrm{Bs}\,}

\newcommand{\Sing}{\mathrm{Sing}\,}

\newcommand{\SigmaA}{\Sigma^{13}_{\mA}}
\numberwithin{equation}{section}

\makeatother

\providecommand{\claimname}{Claim}
\providecommand{\definitionname}{Definition}
\providecommand{\examplename}{Example}
\providecommand{\propositionname}{Proposition}
\providecommand{\theoremname}{Theorem}

\begin{document}
\begin{abstract}
In our previous research, we constructed the affine varieties $\Sigma_{\mA}^{13}$
and $\Pi_{\mA}^{14}$ whose partial projectivizations admit $\mP^{2}\times\mP^{2}$-fibrations
with relative Picard number one. In this paper, we produce prime quasi-smooth
$\mQ$-Fano 3-folds which are anticanonically embedded of codimension
four and belong to 23 (resp.~8) classes in the Graded Ring Database
\cite{GRDB}, as weighted complete intersections in weighted projectivizations
of $\Sigma_{\mA}^{13}$ (resp.~$\Pi_{\mA}^{14}$ or its cone). We
also show that a general member of the anticanonical linear system
of a general prime $\mQ$-Fano $3$-fold constructed in this way is
a quasi-smooth $K3$ surface with at worst Du Val singularities.
\end{abstract}

\maketitle
\markboth{Constructing $\mQ$-Fano $3$-folds}{Hiromichi Takagi}
{\small{}{}}{\small\par}

\section{\textbf{Introduction}}

\subsection{Background}

In this paper, we work over $\mC$, the field of complex numbers.

A complex projective variety is called a \textit{$\mQ$-Fano variety}
if it is a normal variety with only terminal singularities and its
anticanonical divisor is ample. The classification of $\mQ$-Fano
3-folds is one of the central problems in Mori theory for projective
3-folds. In the Graded Ring Database \cite{GRDB}, a huge table of
candidates of $\mQ$-Fano 3-folds is provided. Substantial effort
has been made to show the existence of $\mQ$-Fano 3-folds for such
candidates. 

Hereafter we consider a special class of $\mQ$-Fano 3-folds, called
\textit{prime} $\mQ$-Fano 3-folds, where a $\mQ$-Fano $3$-fold
is said to be \textit{prime} if its anticanonical divisor generates
the group of numerical equivalence classes of $\mQ$-Cartier Weil
divisors on it. The study of prime $\mQ$-Fano 3-folds has a fruitful
history since Gino Fano succeeded in classifying the smooth ones \cite{Fa}. 

To proceed with the classification, we adopt the viewpoint of anticanonical
embeddings into weighted projective spaces. We say that a (not necessarily
prime) $\mQ$-Fano $3$-fold $X$ is \textit{anticanonically embedded
of codimension $c$} in a weighted projective space $\mP$ if the
following conditions are satisfied:

\begin{enumerate}[(1)]

\item $c=\dim\mP-\dim X$, 

\item the anticanonical sheaf $\omega_{X}^{[-1]}$ of $X$ coincides
with $\sO_{\mP}(1)|_{X}$, and

\item the homogeneous coordinate ring of $X$ is Gorenstein.

\end{enumerate} In particular, these conditions imply that the homogeneous
coordinate ring of $X$ coincides with its anticanonical ring $\bigoplus_{n\geq0}H^{0}(X,\omega_{X}^{[-n]})$
(cf.$\,$\cite[(5.16) (ii)]{GW}).

For prime $\mQ$-Fano 3-folds $X$ anticanonically embedded in codimension
$\leq2$, the classification has been completed, provided that $X$
is \textit{quasi-smooth}---that is, its affine cone in the affine
space associated with the anticanonical graded ring is smooth outside
the vertex (\cite{Fl}, \cite{CCC}). Examples of prime $\mQ$-Fano
3-folds anticanonically embedded of codimension $=3$ were constructed
systematically in \cite{Al} using $5\times5$ skew-symmetric matrices. 

\subsection{Prime $\mQ$-Fano $3$-folds anticanonically embedded of codimension
$4$}

As for prime $\mathbb{Q}$-Fano 3-folds anticanonically embedded of
codimension $4$, there are 145 classes of numerical candidates of
them listed in \cite{GRDB}. We exclude the two classes No.$\,$29374
and 26989 in the subsequent explanation: the former is the class of
smooth prime Fano 3-folds of genus 6, classified by Gushel\'\ \cite{Gu};
the latter is the class of prime $\mQ$-Fano $3$-folds of genus $5$
with one $\nicefrac{1}{2}(1,1,1)$-singularity, which we classify
in \cite{Tak1} and revisit in \cite{Tak4}. Several examples of such
prime $\mathbb{Q}$-Fano 3-folds belonging to the remaining 143 classes
are constructed in \cite{CD} as weighted complete intersections in
weighted projectivizations of the $C_{2}$- or $G_{2}^{(4)}$-cluster
variety. We refer informally to this type of construction as \textit{construction
via key varieties}, where the key variety in their setting is the
cluster variety.

In \cite{Tak5}, we introduce the affine variety $\mathscr{H}{}_{\mA}^{13}$,
constructed from the quadratic Jordan algebra of a cubic form, and
show that it contains the $C_{2}$-cluster variety. From $\mathscr{H}{}_{\mA}^{13}$,
we construct more general examples than those in {[}CD{]}, including
prime $\mathbb{Q}$-Fano 3-folds of No.$\,$20652 (No.$\,$5.4 in
\cite{Tak1}), which arise as weighted complete intersections in a
particular weighted projectivization of $\mathscr{H}{}_{{\mathbb{A}}}^{13}$
and are not obtained from the $C_{2}$-cluster variety.

Similarly, we construct the affine variety $\mathfrak{U}_{\mathbb{A}}^{14}$ 
based on the Freudenthal triple system in \cite{Tak7}, which turns out to contain
the $G_{2}^{(4)}$-cluster variety. This again leads to more general
examples than those in \cite{CD}, including prime $\mathbb{Q}$-Fano
3-folds of No.$\,$20544, which are not obtained from the $G_{2}^{(4)}$-cluster
variety.

In \cite{Tak5} (resp.~\cite{Tak7}), we show that a partial projectivization
of $\mathscr{H}_{\mA}^{13}$ (resp.~$\mathfrak{U}_{\mA}^{14}$) admits
a $\mathbb{P}^{2}\times\mathbb{P}^{2}$ (resp.~$\mathbb{P}^{1}\times\mathbb{P}^{1}\times\mathbb{P}^{1}$)-fibration
with relative Picard number one. The above results imply that various
weighted projectivizations of $\mathscr{H}_{\mA}^{13}$ and $\mathfrak{U}_{\mA}^{14}$
produce prime $\mathbb{Q}$-Fano 3-folds of distinct topological types
belonging to the same class in the \cite{GRDB}. This phenomenon can
be regarded as a generalization of that observed in \cite{Tak1} for
prime $\mathbb{Q}$-Fano 3-folds with only $\nicefrac{1}{2}(1,1,1)$-singularities.
One may also view the fact that these various weighted projectivizations
of $\mathscr{H}_{\mA}^{13}$ and $\mathfrak{U}_{\mA}^{14}$ yield
$\mQ$-Fano 3-folds with Picard number one as a consequence of the
property that the above-mentioned partial projectivizations of $\mathscr{H}_{\mA}^{13}$
and $\mathfrak{U}_{\mA}^{14}$ admit the fibrations whose general
fibers have Picard number at least two, while the relative Picard
numbers of the fibrations are one.

Finally in this subsection, we mention the important work \cite{bkr},
in which $\mQ$-Fano 3-folds anticanonically embedded of codimension
$4$ are constructed for 116 out of the 143 classes, although it is
not verified whether they have Picard number one ({[}ibid., the last
para.$\,$of Subsec.$\,$3.4{]}). These varieties are constructed
via the so-called Type I unprojection. For each class, they produce
examples of two distinct topological types, corresponding to the so-called
Tom and Jerry types of unprojection. We believe that these examples
are actually prime.

\subsection{Main result}

In our papers \cite{Tak9,Tak6}, we constructed the affine varieties
$\Sigma_{\mA}^{13}$ and $\Pi_{\mA}^{14}$ respectively, and showed
that their partial projectivizations have $\mP^{2}\times\mP^{2}$-fibrations
with relative Picard number one. 

The purpose of this paper is to construct prime $\mQ$-Fano 3-folds
anticanonically embedded of codimension $4$ belonging to $23$ (resp.~8)
classes in \cite{GRDB}, as weighted complete intersections of weighted
projectivizations of $\Sigma_{\mA}^{13}$ (resp.~$\Pi_{\mA}^{14}$
or its cone). 

To state the main result of this paper, we prepare one convention:
\begin{defn}
\label{def:numdata}By a \textit{numerical data} for a prime $\mQ$-Fano
3-fold $X$ in \cite{GRDB}, we mean the triplet provided in ibid.,
consisting of

\begin{enumerate}[$(a)$]

\item the Hilbert numerator of $X$ with respect to the anticanonical
sheaf $\omega_{X}^{[-1]}$, 

\item the basket of singularities of $X$, and 

\item the weights of the coordinates of the ambient weighted projective
space of $X$.

\end{enumerate}

We remark that the definition of a numerical data here is slightly
different from the one given in \cite[Subsec.\,3.1]{bkr}.
\end{defn}
Note that, from (a) and (c), we may compute 
\begin{quotation}
(d) the Hilbert series $\sum_{n=0}^{\infty}\dim H^{0}(X,\omega_{X}^{[-n]})t^{n}$
of $X$ with respect to $\omega_{X}^{[-1]}$. 
\end{quotation}
Using the data (a)--(c), together with the Riemann-Roch theorem and
the Kawamata-Viehweg vanishing theorem, we may compute important invariants
of $X$; for example, the volume of $\omega_{X}^{[-1]}$ and the genus
$g(X):=h^{0}(\omega_{X}^{[-1]})-2$ can be derived from (b) and (d)
(cf.\cite[(4.6.2)]{abr}). 
\begin{thm}
\label{thm:main1} The following assertions hold$:$ 

\begin{enumerate}[$(1)$]

\item For each of the $23$ (resp.~$8$) numerical data with the
numbers as in Table~\ref{tab:keyvarwt} (resp.~Table~\ref{tab:Weights-of-coordinatesPi})
in Section \ref{sec:Tables}, there exists a quasi-smooth prime $\mQ$-Fano
$3$-fold $X$, anticanonically embedded of codimension $4$, constructed
as a weighted complete intersection in the weighted projectivization
$\Sigma_{\mP}^{12}$ of $\Sigma_{\mA}^{13}$ as in Table~\ref{Table 2}
(resp.~$\Pi_{\mP}^{13}$ of $\Pi_{\mA}^{14}$ or $\Pi_{\mP}^{14}$
of the cone of $\Pi_{\mA}^{14}$ as in Table~\ref{Table3}). The
corresponding coordinate weights are given in Table~\ref{tab:keyvarwt}
and Table~\ref{tab:Weights-of-coordinatesPi}, respectively.

\item A general member $T$ of  $|{-}K_X|$ for a general $X$ as
in (1) is a quasi-smooth $K3$ surface with only Du Val singularities
of type $A$. Moreover, $\Sing T=\Sing X$, and if $X$ has a $1/\alpha\,(\beta,-\beta,1)$-singularity
at a point for some $\alpha,\beta\in\mN$ with $(\alpha,\beta)=1$,
then $T$ has a $1/\alpha\,(\beta,-\beta)$-singularity there. 

\end{enumerate}
\end{thm}
The examples of prime $\mQ$-Fano $3$-folds as in Theorem \ref{thm:main1}
are new. In the paper \cite{Tay}, using Type ${\rm II}_{1}$ unprojection,
the author constructs $\mathbb{Q}$-Fano 3-folds anticanonically embedded
of codimension $4$ with the same numerical data as the one in Theorem
\ref{thm:main1} (1) obtained from $\Sigma_{\mP}^{12}$ except for
the classes No.$\,$393, 642, 644, 4850, 5202, 11004, 16227. Moreover,
the author constructs the $14$-dimensional affine variety and shows
that the examples are obtained in its weighted projectivization as
weighted complete intersections. In our paper \cite{Tak6}, we revisit
the $14$-dimensional variety and call it $\Upsilon_{\mA}^{14}$.
The construction of the affine variety $\Pi_{\mA}^{14}$ was inspired
by \cite[Subsec.\,5.3]{Tay} and the existence of prime $\mathbb{Q}$-Fano
3-folds as in Theorem \ref{thm:main1} (1) obtained from $\Pi_{\mP}^{13}$
or $\Pi_{\mP}^{14}$ was suggested in the examples constructed there
by Type ${\rm II_{2}}$ unprojection except for No.$\,$308. Note
that such an $X$ of No.$\,$308 is shown to be birationally superrigid
by \cite{O}. It is a bonus of our construction via key variety to
confirm the existence of a $\mathbb{Q}$-Fano 3-fold $X$ of No.$\,$308.
We believe that the examples in \cite{Tay} are actually prime. 

Now, examples of prime $\mQ$-Fano 3-folds have been constructed for
141 out of the 143 classes in \cite{GRDB}, using weighted projectivizations
of affine varieties whose partial projectivizations admit $\mP^{2}\times\mP^{2}$-fibrations,
except for the two classes No.$\,$166 and No.$\,$12960. The existence
of prime $\mQ$-Fano 3-folds of No.$\,$166 remains unknown. It is
birationally superrigid if they exist \cite{O}. A quasi-smooth prime
$\mQ$-Fano 3-fold of No.$\,$12960 has six $\nicefrac{1}{2}(1,1,1)$-singularities
and its genus is 1. Recently, in \cite{Tak10}, we construct examples
of such $\mQ$-Fano $3$-folds with two distinct topological types
via other key varieties. Partial projectivization of one of the key
varieties admits a fibration with relative Picard number one whose
general fibers are isomorphic to the $5$-dimensional cone over $\mP^{2}\times\mP^{2}$.

We expect that the affine variety $\mathfrak{U}_{\mA}^{14}$ will
produce more examples of prime $\mQ$-Fano 3-folds anticanonically
embedded of codimension $4$. We plan to study this problem in the
near future.

\vspace{3pt}

The structure of the paper is as follows: In Section \ref{sec:Strategy-to-show},
we explain in detail the strategy for proving Theorem \ref{thm:main1}.
We then prove the theorem in Section \ref{sec:Proof-of-theorem1}
for $\Sigma_{\mP}^{12}$ and in Section \ref{sec:Proof-of-theorem2}
for $\Pi_{\mP}^{13}$ or $\Pi_{\mP}^{14}$. For the reader's convenience,
we review the properties of the affine varieties $\Sigma_{\mA}^{13}$
and $\Pi_{\mA}^{14}$ in Section \ref{sec:Equations-of-the key variety}. 

\vspace{5pt}

As explained in Section \ref{sec:Strategy-to-show}, the analysis
of the singularities of the anticanonical $K3$ surface $T$ in each
case plays a crucial role in the proof of Theorem \ref{thm:main1}.
Although the computations for the analysis are difficult to carry
out by hand, they are in fact straightforward and can be performed
without difficulty using built-in commands in the software systems
\textit{Magma }\cite{BCP} and\textit{ Mathematica} \cite{W}. Part
of the relevant\textit{ Mathematica} code is available in \cite{Tak8}.
The reason why these computations are straightforward lies in construction
via key varieties of the $K3$ surface $T$. Indeed, in the process
of obtaining the $K3$ surface $T$ as a section of the corresponding
higher-dimensional key variety, we encounter intermediate varieties
whose affine cones are complete intersections in some coordinate charts.
Therefore this allows us to verify the singularities using the Jacobian
criterion.

\section{\textbf{Strategy to show Theorem \ref{thm:main1}} \label{sec:Strategy-to-show}}

In this section, we explain in detail our strategy how to verify weighted
complete intersections of key varieties are desired examples of prime
$\mQ$-Fano $3$-folds and $K3$ surfaces. In Sections \ref{sec:Proof-of-theorem1}
and \ref{sec:Proof-of-theorem2}, we carry out the verification following
the strategy explained in this section. 

\subsection{Set-up\label{subsec:Set-up}}

\vspace{3pt}

Fix a class in Tables~\ref{tab:keyvarwt} and~\ref{tab:Weights-of-coordinatesPi}
in Section \ref{sec:Tables}. Hereafter in this section, we denote
by $\mathfrak{K}$ the corresponding weighted projective key variety
and by $\mathfrak{K}_{\mA}$ the affine cone of $\mathfrak{K}$ ,
namely $\mathfrak{K}_{\mA}=\Sigma_{\mA}^{13},\Pi_{\mA}^{14}$ or the
cone over $\Pi_{\mA}^{14}$. In each case, we set 
\begin{equation}
X=\mathfrak{K}\cap(a_{1})^{m_{1}}\cap\dots\cap(a_{k})^{m_{k}}\label{eq:Xgen}
\end{equation}
where $(a_{i})$ for each $i$ means a general weighted hypersurface
of weight $a_{i}$, called a \textit{section} of weight $a_{i}$.
The numbers $a_{1},\dots,a_{k}$, $m_{1},\dots,m_{k}\in\mN$ are assigned
for the class, and may be assumed to satisfy $a_{1}<\cdots<a_{k}$
and $m_{1}+\cdots+m_{k}=\dim\mathfrak{K}-3$.

Let $T$ be the intersection of $X$ and a section of weight one.
We may write 

\[
T=\begin{cases}
\mathfrak{K}\cap(a_{1})^{m_{1}+1}\cap\dots\cap(a_{k})^{m_{k}} & :a_{1}=1,\\
\mathfrak{K}\cap(1)\cap(a_{1})^{m_{1}}\cap\dots\cap(a_{k})^{m_{k}} & :a_{1}>1.
\end{cases}
\]
We write this uniformly as
\begin{equation}
T=\mathfrak{K}\cap(b_{1})^{n_{1}}\cap\dots\cap(b_{l})^{n_{l}}\label{eq:Tgen}
\end{equation}
with $b_{1}<b_{2}<\cdots<b_{l}$. 

\vspace{3pt}

\noindent \textbf{List of notation}:

We use the following notation for a weighted projective variety $V$
considered below (for example, $V=X,T$): 
\begin{itemize}
\item $\mP_{V}:=$ the ambient weighted projective space for $V$. 
\item The $x$-chart of $V$$:=$ $V\cap\{x\not=0\}$ for a coordinate $x$
of $\mP_{V}$.
\item The $x$-point of $V$$:=$the point of $V$ whose only nonzero coordinate
is $x$. 
\item $V_{\mA}:=$ the affine cone of $V$. 
\item $V_{\mA}^{o}:=$ the complement of the origin in $V_{\mA}$. 
\item $\mA_{V}:=$ the ambient affine space for $V_{\mA}$. 
\item $\mA_{V}^{x}:=\{x=1\}\subset\mA_{V}$. 
\item $V_{\mA}^{x}:=V_{\mA}\cap\{x=1\}.$ We also call this the $x$-\textit{chart}
of $V_{\mA}$. Note that $V_{\mA}^{x}$ is a closed subset of $V_{\mA}^{o}$.
\item $|\sO_{\mP_{V}}(i)|_{\mA}:=$ the set of the affine cones of members
of $|\sO_{\mP_{V}}(i)|$. 
\item $\Bs|\sO_{\mP_{V}}(i)|_{\mA}:=$ the intersection of all the members
of $|\sO_{\mP_{V}}(i)|_{\mA}$. 
\end{itemize}

\subsection{Shape of the equations of $X$ and $T$\label{subsec:Shape-of-equations}}

In each case, it is important to observe that the number of the coordinates
of $\mP_{\mathfrak{K}}$ with weight $a_{i}$ is greater than or equal
to $m_{i}$ for any $i$ (we can read off this fact from Tables ~\ref{tab:keyvarwt}
and~\ref{tab:Weights-of-coordinatesPi} in Section \ref{sec:Tables}).
Therefore we have the following claim:
\begin{claim}
\label{claim:ShapeX}For each $X$, we may choose each of the sections
$(a_{i})$ as defined by the equation of the form 
\begin{align}
 & \text{(a coordinate of \ensuremath{\mP_{\mathfrak{K}}} with weight \ensuremath{a_{i})=}}\label{sectioneq}\\
 & \quad\text{(a polynomial with weight \ensuremath{a_{i}} of other coordinates)}.\nonumber 
\end{align}
\end{claim}
Substituting the r.h.s.\,of the equations of the sections of the
form \eqref{sectioneq} for the corresponding variables of the nine
equations of $\mathfrak{K}$, we obtain the equations of $X$ in their
ambient weighted projective spaces $\mP_{X}$.

Moreover, since the number of the coordinates of $\mP_{\mathfrak{K}}$
with weight one is greater than or equal to $n_{1}$ as in \eqref{eq:Tgen},
we similarly obtain the equations of $T$ in their ambient weighted
projective spaces $\mP_{T}$. 

\subsection{Presentation of the equations of $X$ and $T$\label{subsec:Presentation-of-the eq XT}}

Usually, the equations of $X$ as stated in the subsection \ref{subsec:Shape-of-equations}
are complicated due to the presence of many parameters. Fortunately,
to show Theorem \ref{thm:main1}, we do not need the explicit equations
of $X$ in any case. 

In the case that $h^{0}(\sO_{\mP_{X}}(1))=1$, the equations of $T$,
which is relatively simple, suffice to show Theorem \ref{thm:main1}
both for $X$ and $T$. The equations of $T$ in each such case are
given in Subsection \ref{subsec:Caseh0=00003D1} and Section \ref{sec:Proof-of-theorem2}. 

In the case that $h^{0}(\sO_{\mP_{X}}(1))=2$, the full equations
of $T$ are still complicated. To overcome this situation, we take
two sections $T,T'$ of $X$ by general members of $|\sO_{\mP_{X}}(1)|$
and set $C:=T\cap T'$. It turns out that the equations of $C$, which
are now relatively simple, suffice to show Theorem \ref{thm:main1}
both for $X$ and $T$. The equations of $C$ in each such case are
given in Subsection \ref{subsec:Caseh0=00003D2}. 

In the case that $h^{0}(\sO_{\mP_{X}}(1))\geq3$, it turns out that
we need only a small amount of information about the equations of
$X$ and $T$ to show Theorem \ref{thm:main1}.

\subsection{Main checkpoints: Claims (A), (B), (C)}

\noindent For each case, we will verify the following Claims (A)--(C)
for general $X$ and $T$: 

\vspace{5pt}

\noindent \textbf{Claim (A).} \textit{$X$ is a quasi-smooth $3$-fold
and $T$ is a quasi-smooth surface; namely, $X_{\mA}^{o}$ is a smooth
$4$-fold and $T_{\mA}^{o}$ is a smooth $3$-fold. }

\vspace{5pt}

\noindent \textbf{Claim (B).} \textit{It holds that $\Sing X=\Sing T$.
$X$ (resp.$\,T$) have only cyclic quotient singularities assigned
as in Tables~\ref{Table 2} and~\ref{Table3} in Section \ref{sec:Tables}
(resp.$\,$as in Theorem \ref{thm:main1} (2)). Moreover, if $X$
has a $1/\alpha(1,\beta,\alpha-\beta)$-singularity for some coprime
$\alpha$ and $\beta$ with $0<\beta<\alpha$, then we may choose
$\beta$ is the residue of the weight of a coordinate by $\alpha$.
A similar assertion holds for $T$.}

\vspace{5pt}

\noindent \textbf{Claim (C).} \textit{We set 
\[
\mathsf{b}:=\begin{cases}
p_{1} & :\mathfrak{K}=\Sigma_{\mP}^{12}\\
p_{3}^{2}+t_{1}p_{4}^{2}+p_{4}u_{2} & :\mathfrak{K}=\Pi_{\mP}^{13}\ \text{or}\ \Pi_{\mP}^{14}.
\end{cases}
\]
It holds that $X\cap\{\mathsf{b}=0\}$ is a prime divisor.} 

\textit{We refer to Proposition \ref{prop:Pic1} for the meaning of
$\mathsf{b}$. We sometimes call $X\cap\{\mathsf{b}=0\}$ the }boundary.\textit{ }

\vspace{5pt}

\subsection{Claims (A)--(C) imply Theorem \ref{thm:main1}\label{subsec:Claims-(A)=002013(C)-imply MainThm}}

By Claims (A)--(C), we can finish the proof of Theorem \ref{thm:main1}
as follows: 

\vspace{3pt}

\noindent \textbf{Proof of (1).} We show that $X$ is a prime $\mQ$-Fano
3-fold. Claim (A) implies that $X$ is irreducible since $\mathfrak{K}$
is irreducible and normal by Proposition \ref{prop:Sigmamain} and
$X$ is a weighted complete intersection in $\mathfrak{K}$ cut out
by ample divisors. By Claim (B), $X$ has only terminal singularities.
By Claims (A), (C) and Proposition \ref{prop:Pic1} (with the remark
following it), $X$ has Picard number one. By Proposition \ref{prop:candiv},
we can immediately compute the canonical sheaf of $X$, and then we
see that $\omega_{X}=\sO_{X}(-1)$, which is primitive by Proposition
\ref{prop:Pic1}. Thus $X$ is a prime $\mQ$-Fano 3-fold.

We show that $X$ is anticanonically embedded of codimension $4$.
Indeed, $X$ is of codimension $4$ in the ambient weighted projective
space since so is $\mathfrak{K}$. The condition that $\omega_{X}^{[-1]}=\sO_{X}(1)$
has been already verified above. The homogeneous coordinate ring of
$X$ is Gorenstein since so is that of $\mathfrak{K}$ by Proposition
\ref{prop:Sigmamain} (1) and $X$ is general (cf.$\,$\cite[Prop.\,1.3]{R}).
Therefore $X$ is anticanonically embedded of codimension $4$. 

Finally, we verify the numerical data (a)--(c) as in Definition \ref{def:numdata}
coincides with those given in \cite{GRDB} for each $X$. Let $X$
be obtained from $\mathfrak{K}=\Sigma_{\mP}^{12}$ (resp.~$\mathfrak{K}=\Pi_{\mP}^{13}\ \text{or}\ \Pi_{\mP}^{14}$).
The data (c) can be read off from the weights assigned to the coordinates
of $\mathfrak{K}$ and the way $X$ is cut out from it. The data (b)
can be read off from Claim (B). As for the data (a), it can be read
as follows. First, from the free resolution \cite[(2.17)]{Tak9} (resp.$\,$\cite[(5.7)]{Tak6})
of the homogeneous ideal of $\mathfrak{K}=\Sigma_{\mP}^{12}$ (resp.$\,$$\mathfrak{K}=\Pi_{\mP}^{13}\ \text{or}\ \Pi_{\mP}^{14}$),
we obtain the Hilbert series of its homogeneous coordinate ring, and
hence the Hilbert series of the homogeneous coordinate ring of $X$.
Note that the homogeneous coordinate ring of $X$ coincides with its
anticanonical ring since the homogeneous coordinate ring of $X$ has
been proved to be Gorenstein above (cf.$\,$$\,$\cite[(5.16) (ii)]{GW}).
Therefore the data (a) is obtained from this together with data (c).
Then we can verify that the data (a)--(c) obtained in this way indeed
agree with those listed in {[}GRDB{]}.

Thus we obtain Theorem \ref{thm:main1} (1). 

\noindent \textbf{Proof of (2).} Since $T\in |{-}K_X|$, we have
$K_{T}\sim0$ and $h^{1}(\sO_{T})=0$. Thus, by Claims (A) and (B),
$T$ is a desired quasi-smooth $K3$ surface, which shows Theorem
\ref{thm:main1} (2).

\subsection{Singularities and base loci}

In this subsection, we prove the following claim, which is important
for several considerations in the sequel:
\begin{claim}
\noindent \label{claim:SingBase}It holds that
\begin{equation}
\Sing X_{\mA}^{o}\subset\bigcup_{i}\Bs|\sO_{\mP_{\mathfrak{K}}}(a_{i})|_{\mA},\quad\Sing T_{\mA}^{o}\subset\bigcup_{i}\Bs|\sO_{\mP_{\mathfrak{K}}}(b_{i})|_{\mA}.\label{eq:SingXT}
\end{equation}
\end{claim}
\begin{proof}
\noindent By the Bertini singularity theorem (cf.\cite{Kl}), we see
that 

\noindent 
\begin{align*}
\Sing X_{\mA}^{o} & \subset\left(\bigcup_{i}\Bs\left|\sO_{\mP_{\mathfrak{K}}}(a_{i})\right|_{\mA}\right)\cup\Sing\mathfrak{K}_{\mA}^{o},\\
\Sing T_{\mA}^{o} & \subset\left(\bigcup_{i}\Bs\left|\sO_{\mP_{\mathfrak{K}}}(b_{i})\right|_{\mA}\right)\cup\Sing\mathfrak{K}_{\mA}^{o}.
\end{align*}
 Note that the dimension of the singular locus of $\mathfrak{K}_{\mA}$
is less than the codimension of $X$ in $\mathfrak{K}$ by \cite[Prop. 2.14]{Tak9}
for $\mathfrak{K}=\Sigma_{\mP}^{12}$ and \cite[Prop. 6.1]{Tak6}
for $\mathfrak{K}=\Pi_{\mP}^{13}\ \text{or}\ \Pi_{\mP}^{14}$. Therefore
we see that $(\Sing\mathfrak{K}_{\mA}^{o})\setminus(\bigcup_{i}\Bs|\sO_{\mP_{\mathfrak{K}}}(a_{i})|_{\mA})$
is disjoint from $X_{\mA}^{o}$, and $(\Sing\mathfrak{K_{\mA}^{o}})\setminus(\bigcup_{i}\Bs|\sO_{\mP_{\mathfrak{K}}}(b_{i})|_{\mA})$
is disjoint from $T_{\mA}^{o}$ for general $X$ and $T$. Therefore
the relations in \eqref{eq:SingXT} hold. 
\end{proof}

\subsection{Reduction of Claims (A) and (B) to the case of $T$\label{subsec:Reduction-of-Claims}}

For Claims (A) and (B), we show the following reduction step:
\begin{claim}
\label{cla:RedAB} Claim (A) for $T$ implies Claim (A) for $X$.
The equality $\Sing X=\Sing T$ in Claim (B) also holds. As for the
rest of Claim (B), it suffices to verify the corresponding assertions
for $T$.
\end{claim}
\begin{proof}
\noindent \textbf{Claim (A):} It suffices to show $\Sing X_{\mA}^{o}\subset\Sing T_{\mA}^{o}$.
Since $T_{\mA}^{o}$ is a Cartier divisor of $X_{\mA}^{o}$, we have
\begin{equation}
(\Sing X_{\mA}^{o})\cap T_{\mA}^{o}\subset\Sing T_{\mA}^{o}.\label{eq:SingXSingT}
\end{equation}
Since there exists a coordinate of $\mA_{\mathfrak{K}}$ with weight
one, $\Bs|\sO_{\mP_{\mathfrak{K}}}(a_{i})|_{\mA}$ is contained in
$\Bs|\sO_{\mP_{\mathfrak{K}}}(1)|_{\mA}$ for any $i$, hence, by
Claim \ref{claim:SingBase}, $\Sing X_{\mA}^{o}$ is contained in
$\Bs|\sO_{\mP_{\mathfrak{K}}}(1)|_{\mA}\cap X_{\mA}^{o}$. Moreover,
since $T_{\mA}^{o}$ is cut from $X_{\mA}^{o}$ by a general section
of weight one, $\Bs|\sO_{\mP_{\mathfrak{K}}}(1)|_{\mA}\cap X_{\mA}^{o}$
is contained in $T_{\mA}^{o}$. Therefore $\Sing X_{\mA}^{o}$ is
contained in $T_{\mA}^{o}$, and hence, by \eqref{eq:SingXSingT},
we obtain $\Sing X_{\mA}^{o}\subset\Sing T_{\mA}^{o}$ as desired. 

\vspace{3pt}

\noindent \textbf{Claim (B):} By Claim (A), $X$ and $T$ has singularities
along the quotient of the nonfree locus of the $\mC^{*}$-action on
$X_{\mA}^{o}$ and $T_{\mA}^{o}$ respectively. Moreover, the nonfree
locus of the $\mC^{*}$-action on $X_{\mA}^{o}$ is contained in ${\rm \Bs|\sO_{\mP_{\mathfrak{K}}}(1)|_{\mA}}$,
hence in $T_{\mA}^{o}$. Therefore it holds that $\Sing X={\rm \Sing}T$.
If $T$ has the singularity of the desired type $1/\alpha\,(\beta,\alpha-\beta)$
at a singular point ${\sf t}$ and $\beta$ is the residue of the
weight of a coordinate by $\alpha$, then $X$ has the desired $1/\alpha(1,\beta,\alpha-\beta)$-singularity
at ${\sf t}$ since $T=X\cap(1)$. This shows that the remaining part
of Claim (B) follows from the corresponding assertions for $T$.
\end{proof}
\vspace{10pt}

The proofs of Claims (A) and (B) for $T$ are outlined at the beginning
of Subsections \ref{subsec:Caseh0=00003D1}, \ref{subsec:Caseh0=00003D2},
and \ref{subsec:Caseh0ge3}, and Section \ref{sec:Proof-of-theorem2}.
Since the essential structure of the argument remains the same for
each case, we provide complete proofs, accompanied by the full \textit{Mathematica}
codes \cite{Tak8}, only for the cases No.$\,$360, 1185 and 577.
For the other cases, only the results of the computations are recorded
as data.

\vspace{3pt}

In what follows, we describe some techniques employed in the proofs
of Claims (A) and (B).

\subsection{Linear Part Computation (LPC)\label{subsec:Linear-Part-Computation}}

\noindent Taking into account the weights of the coordinates of the
ambient weighted projective space $\mP_{T}$ of $T$, it is possible
to determine in advance where $T$ should have singularities. The
method LPC, explained below, is applied mainly at such points to determine
all the singularities of $T$ except for $1/2(1,1)$-singularities,
which are treated separately (see Subsection \ref{subsec:1/2-singularities}
for an explanation). LPC also serves to demonstrate the smoothness
of $T_{\mA}^{o}$ without using the Jacobian criterion, hence it is
also used to show that $T$ is smooth at some point. 

\subsubsection{\textbf{LPC for $T$}}

We first select a point ${\sf t}\in T$ at which we intend to apply
the method LPC, to be described below. We choose a coordinate $x$
of $\mP_{T}$ which is not zero for the point $\sf{t}$, and one point
${\sf t'}\in T_{\mA}^{x}$ corresponding to ${\sf t}$. We localize
the equations of $T$ at ${\sf t'}$ setting $x=1$ and transforming
the coordinates of $\mA_{T}^{x}$ to make the point ${\sf t}'$ the
origin. For example, in the case of the $u$-point, we have only to
set $u=1$. We compute the linear parts of the localized equations
of $T$ and show that they span a four-dimensional subspace in the
cotangent space of ${\mA}_{T}^{x}$ at ${\sf t}'$ (the name LPC comes
from this computation). Then we check that there are two local coordinates
of $T_{\mA}^{x}$ at ${\sf t}'$. This implies that $T$ is quasi-smooth
at ${\sf t}$. Let $\alpha$ be the order of the stabilizer group
at ${\sf t}'$ of the group action induced on $\mA_{T}^{x}$ from
$\mP_{T}$. In each case, we see that the weights of the two local
coordinates with respect to the $\mZ/\alpha\mZ$-action are $\beta$
and $\alpha-\beta$ with some $\beta\in\mN$ coprime to $\alpha$,
and $\beta$ is the residue of the weight of a coordinate by $\alpha$.
Thus $T$ has a $1/\alpha(\beta,\alpha-\beta)$-singularity at ${\sf t}$
which is the desired type in each case.

While the nonsingularity of $T_{\mA}^{x}$ at ${\sf t'}$ could be
verified by the Jacobian criterion, as will be subsequently explained
at the beginning of Subsection \ref{subsec:Caseh0=00003D1}, the advantage
of LPC is that, as we have seen above, it also determines the singularity
type of $T$ at $\sf{t}$.

\vspace{5pt}

We present the following example illustrating the application of LPC.
In this example, we refer to the data of No.$\,$360, which will be
explicitly given in Subsection \ref{subsec:Caseh0=00003D1}; the reader
is encouraged to consult that data while reading the discussion. 

\begin{example}
\label{exa:LPCNo.360.}We explain how to determine the singularities
of $T$ of No.$\,$360 on the $p_{2}$-chart by LPC. Since $w(p_{2})=8$,
the cyclic group $\mZ_{8}$ acts on $\mA_{T}^{p_{2}}$. By the data
\eqref{eq:360Tamb}, we see that the action on $\mA_{T}^{p_{2}}$
is not free along $\{t_{1}=p_{1}=r=p_{3}=0\}$. By the explicit equations
of $T$ (see data \textbf{Sections for $T$} given in Subsection \ref{subsec:Caseh0=00003D1}
and the equations of $\Sigma_{\mA}^{13}$ as in Definition \ref{defn:Key})
we see that this intersects $T_{\mA}^{p_{2}}$ along $\{-a_{11}a_{3}^{2}t_{2}^{4}+b_{0}t_{2}^{2}+a_{0}=0\}$
in the $t_{2}$-line, which map to two points on $T$ by the $\mZ_{8}$-quotient.
We adopt LPC for these two points. For this, we set $t_{2}=1$ instead
of setting $p_{2}=1$. In $\mA_{T}^{t_{2}}$, the inverse image of
these two points on $T$ are $T_{\mA}^{t_{2}}\cap\{-a_{11}a_{3}^{2}+b_{0}p_{2}+a_{0}p_{2}^{2}=0\}$
in the $p_{2}$-line, which consists of two points. Choose one $A\in\mC$
such that $A^{2}=b_{0}^{2}+4a_{0}a_{11}a_{3}^{2}$. Then the $p_{2}$-coordinates
of these two points in $T_{\mA}^{t_{2}}$ are $\frac{-b_{0}\pm A}{2a_{0}}$.
By the corresponding coordinate change $p_{2}=p'_{2}+\frac{-b_{0}+A}{2a_{0}}$,
or $p_{2}=p'_{2}+\frac{-b_{0}-A}{2a_{0}}$ for each of two points,
it becomes the origin of the $\mA_{T}^{t_{2}}$. Then we can apply
LPC. Consequently, we see that $T$ has $1/4(1,3)$-singularities
at these two points. Details can be found in the \textit{Mathematica}
code for No.$\,$360 \cite{Tak8}.
\end{example}

\subsubsection{\textbf{LPC for $C$\label{subsec:LPC-for C}}}

\noindent The method LPC works also for $C$, which is defined as
in Subsection \ref{subsec:Presentation-of-the eq XT} when $h^{0}(\sO_{\mP_{X}}(1))=2$.
We fix a point ${\sf t}\in C$ at which we intend to apply LPC. We
choose a coordinate $x$ of $\mP_{C}$ which is not zero for $\sf{t}$,
and one point ${\sf t'}\in C_{\mA}^{x}$ corresponding to ${\sf t}$.
Exactly in the same way as LPC for $T$, we see that there is one
local coordinate of $C_{\mA}^{x}$ at ${\sf t}'$. Let $\alpha$ be
the order of the stabilizer group at ${\sf t}'$ of the group action
induced on $\mA_{C}^{x}$ from $\mP_{C}$. In each case, we observe
that the weight of the local coordinate with respect to the $\mZ/\alpha\mZ$-action
is $\alpha-1$ and $\alpha-1$ is the residue of the weight of a coordinate
by $\alpha$. Then, $T$ has a $1/\alpha\,(1,\alpha-1)$-singularity
at ${\sf t}$.
\begin{example}[\textbf{LPC at the $u$-point}]
\noindent \textbf{\label{exa:LPC at u-pt}}It turns out that $T$
contains the $u$-point in all the cases except for No.$\,$360 and
No.$\,$1218. For them, we adopt LPC at the $u$-point to show $T$
has the desired type of singularity there. To perform LPC at the $u$-point,
we do not need the full equations of the sections defining $T$ or
$C$. We explain this reason for $T$. As for the equations of the
sections defining $T$, we have only to present linear parts of them.
Indeed, the inverse image of the $u$-point is the origin of $\mA_{T}^{u}$
if we set $u=1$. We observe that, in each case, $b_{i}$ for any
$i\geq1$ is less than $w(u)$, hence $u$ does not appear in the
equations of the sections as in \refeq{sectioneq}. Therefore, by
setting $u=1$, these equations of the sections do not change. The
higher degree parts of them do not contribute to the linear parts
of the equations of $T$ with $u=1$, hence we may ignore them when
we adopt LPC. This observation is particularly useful in the case
that $h^{0}(\sO_{X}(1))\geq3$.

As an explicit application example, we consider No.$\,$1185. Note
that the inverse image of the $u$-point is the origin $o$ of $\mA_{T}^{u}$,
and the coordinates of $\mA_{T}^{u}$ are $s_{13},s_{12},p_{1},p_{2},t_{2},p_{3}$
by the data \eqref{eq:1185X}. We may easily compute that the linear
parts of the equations of $T$ with $u=1$ generates the subspace
spanned by $s_{13},p_{1},p_{2},p_{3}$ of the cotangent space of $\mA_{T}^{u}$
at $o$. Thus, as a local coordinate of $T_{\mA}^{u}$ at the point
$o$, we may take $s_{12},t_{2}$. Noting that $w(s_{12})=3$, and
$w(t_{2})=5$, we see that $T$ has a $1/8(3,5)$-singularity at the
$u$-point. Details can be found in the \textit{Mathematica} code
for No.$\,$1185.
\end{example}
\vspace{5pt}

\subsection{$1/2(1,1)$-singularities\label{subsec:1/2-singularities}}

\noindent It turns out that the surface $T$ has $1/2(1,1)$-singularities
in many cases. Usually, the coordinates of $1/2(1,1)$-singularities
are complicated and hence LPC is not appropriate for them. Nevertheless,
the Jacobian criterion is sufficient for them since we can immediately
conclude that a point is a $1/2(1,1)$-singularity once we know it
is obtained locally by the $\mZ_{2}$-quotient of a smooth surface
with an isolated fixed point. An explicit example is calculated in
\textit{Mathematica} code for No.$\,$1185 \cite{Tak8}. 

\subsection{Strategy to prove Claim (C)\label{subsec:(A) and (B) =00003D>(C)}}

We outline the proof of Claim (C), assuming that Claims (A) and (B)
have been established.

As we have already seen in Subsection \ref{subsec:Claims-(A)=002013(C)-imply MainThm}
(Proof of (1)), Claim (A) and (B) implies that $X$ is an irreducible
normal $3$-fold. Therefore, to prove $X\cap\{\mathsf{b}=0\}$ is
irreducible, it suffices to show that it is normal since $X\cap\{\mathsf{b}=0\}$
is an ample divisor of $X$. Since $X$ has only terminal singularities
by Claim (B), $X\cap\{\mathsf{b}=0\}$ is Cohen-Macaulay (cf.~\cite[Cor.\,5.25]{KM}).
Therefore it suffices to show that 
\begin{equation}
\dim\Sing(X\cap\{\mathsf{b}=0\})\leq0.\label{eq:0dimT1}
\end{equation}
Similarly to the proof of Claim \ref{cla:RedAB}, we now reduce \eqref{eq:0dimT1}
to the following claim on $T$:

\begin{equation}
\text{\ensuremath{\dim\Sing(T\cap\{\mathsf{b}=0\})\leq0}}.\label{0dimT}
\end{equation}
Indeed, note that $\mathfrak{K}\cap\{\mathsf{b}=0\}$ is normal by
Proposition \ref{prop:Pic1}. Therefore, similarly to the proof of
Claim \ref{claim:SingBase}, we see that \eqref{eq:0dimT1} holds
outside $\bigcup_{i}\Bs\left|\sO_{\mP_{\mathfrak{K}}}(a_{i})\right|$.
Since there exists a coordinate of $\mA_{\mathfrak{K}}$ with weight
one, $\Bs|\sO_{\mP_{\mathfrak{K}}}(a_{i})|$ is contained in $\Bs|\sO_{\mP_{\mathfrak{K}}}(1)|$
for any $i$, hence, we have only to check \eqref{eq:0dimT1} along
$T$. Since $T$ is a Cartier divisor on $X$ outside a finite set
of points by Claim (B), we have 
\[
\left({\rm \Sing(X\cap\{\mathsf{b}=0\})}\right)\cap T\subset\Sing(T\cap\{\mathsf{b}=0\})
\]
outside a finite set of points. Consequently, it suffices to show
\eqref{0dimT} for \eqref{eq:0dimT1}. 

\vspace{10pt}

The proof of \eqref{0dimT} can be carried out in parallel with that
of Claim (A). It is also outlined at the beginning of Subsections
\ref{subsec:Caseh0=00003D1}, \ref{subsec:Caseh0=00003D2}, and \ref{subsec:Caseh0ge3},
and Section \ref{sec:Proof-of-theorem2}, and complete proofs are
provided only for the cases No.$\,$360, 1185 and 577 in \textit{Mathematica}
code \cite{Tak8}. 

\section{\textbf{Proof of Theorem \ref{thm:main1} for Table~\ref{Table 2}\label{sec:Proof-of-theorem1}}}

\subsection{Case $h^{0}(\sO_{\mP_{X}}(1))=1$\label{subsec:Caseh0=00003D1}}

To establish Claims (A) and (B) for each class of $T$, we proceed
by dividing our analysis into the $p_{1}$-chart, the $p_{2}$-chart,
and the locus $\{p_{1}=p_{2}=0\}|_{T}$. In what follows, we outline
the analysis for each of these cases.

\noindent \vspace{3pt}

\noindent \noindent\textbf{On the $p_{1}$- and $p_{2}$-charts:}
In each case, the smoothness of the $p_{1}$- and $p_{2}$-charts
of $T_{\mA}$ is verified successively in this order. This order of
analysis allows us to assume $p_{1}=0$ when studying the $p_{2}$-chart,
which reduces the amount of computation considerably. Moreover, this
approach aligns conveniently with the analysis of singularities of
the boundary $\{p_{1}=0\}|_{T}$.

Working out the explicit equations in each case, we find that the
$p_{1}$- and $p_{2}$-charts of $T_{\mA}$ are complete intersections
of codimension at most two (some cases being the affine space itself)
except for the $p_{1}$-charts of $T_{\mA}$ for No.$\,$1185 and
No.$\,$1218, and the $p_{1}$- and $p_{2}$-charts of $T_{\mA}$
for No.$\,$1413. For instance, the defining ideal of the $p_{1}$-chart
of $T_{\mA}$ for No.$\,$393 is generated by the equations $F_{1}|_{T}$
and $F_{4}|_{T}$ (the polynomials $F_{i}\,(1\leq i\leq9)$ are defined
as in Subsection \ref{subsec:The-key-varietySigma} below). This was
verified by directly checking in \textit{Magma} \cite{BCP} that the
remaining equations lie in the ideal generated by $F_{1}|_{T}$ and
$F_{4}|_{T}$. 

When the $p_{1}$- or $p_{2}$-chart of $T_{\mA}$ is a complete intersection,
the Jacobian criterion may be used to determine its smoothness in
principle. Indeed, in the $p_{1}$-chart of $T_{\mA}$ for No.$\,$393,
we consider the ideal generated by $F_{1}|_{T}$ and $F_{4}|_{T}$,
and the Jacobian ideal. Computing its elimination ideal in \textit{Magma}
\cite{BCP}, we find that it is nontrivial as for the parameters;
this implies that, for general choices of the parameters, $T_{\mA}$
is nonsingular on the $p_{1}$-chart. Nonetheless, even when a chart
is a complete intersection, elimination computations may often remain
difficult for it. Below, we present two strategies to overcome such
cases.

\vspace{3pt}

\noindent\textit{First computational strategy:} We consider No.$\,$574
for instance. Since $w(p_{1})=5$, the $p_{1}$-chart of $T_{\mA}$
does not intersect ${\rm Bs}|\sO_{\mP_{T}}(5)|$. By the Bertini theorem,
if the variety $U_{\mA}$, obtained by removing the section of weight
$5$ from the equations of $T_{\mA}$, is nonsingular on the $p_{1}$-chart,
then so is $T_{\mA}$. Working with $U_{\mA}$ facilitates explicit
computations, not only because it involves fewer parameters, but also
because the $p_{1}$-chart of $U_{\mA}$ is a hypersurface, in contrast
to that $T_{\mA}$ is a codimension two complete intersection. In
fact, the smoothness of the $p_{1}$-chart of $U_{\mA}$ can be readily
checked using \textit{Mathematica} \cite{W}.

This approach is also effective for No.$\,$1218. Although the $p_{1}$-chart
of $T_{\mA}$ for this class is not a complete intersection, the variety
$U_{\mA}$, obtained by removing the two sections of weight $4$ from
the equations of $T_{\mA}$, has the $p_{1}$-chart that is a hypersurface
and whose smoothness is easily verified. Since $w(p_{1})=4$, the
Bertini theorem again ensures that $T_{\mA}$ is nonsingular on the
$p_{1}$-chart.

\vspace{3pt}

\noindent\textit{Second computational strategy: }This method is illustrated
by No.$\,$360. Removing the section of weight $8$ from the equations
of $T_{\mA}$, we obtain the variety $U_{\mA}$ whose $p_{1}$-chart
is a hypersurface (in contrast to that $T_{\mA}$ is a codimension-two
complete intersection on the $p_{1}$-chart), and whose smoothness
is easily checked in \textit{Mathematica} \cite{W}. Up to this point,
the procedure is similar to the first strategy. However, the key difference
lies in the weight: since $w(p_{1})=7\not=8$, one cannot immediately
conclude that $T_{\mA}$ is nonsingular on the $p_{1}$-chart. Nevertheless,
since we only need to check singularities along the base locus ${\rm Bs}|\sO_{\mP_{T}}(8)|_{\mA}$,
the computational complexity is significantly reduced. In practice,
we perform this check on $T$. It turns out that ${\rm Bs}|\sO_{\mP_{T}}(8)|\bigcap T$
consists of the $p_{1}$-point and the $p_{4}$-point, and singularities
at these point can be determined by conducting LPC. For details, see
the data for No.$\,$360 below and the provided \textit{Mathematica}
code \cite{Tak8}.

\noindent \vspace{3pt}

\noindent \noindent\textbf{The locus $\{p_{1}=p_{2}=0\}_{|T}$}:
In each case, we can easily compute the locus $\{p_{1}=p_{2}=0\}_{|T}$.
Then we observe that $\{p_{1}=p_{2}=0\}_{|T}$ consists of finite
number of points. To show Claims (A) and (B) on $\{p_{1}=p_{2}=0\}|_{T}$,
LPC at these points is sufficient for all the cases. For No.$\,$569,
1185, 1186, and 1218, however, the LPC computations are much more
involved. Therefore, for these cases, we adopt a different approach.
Below, we explain this method only for one representative example,
No.$\,$1185. The corresponding \textit{Mathematica} code is also
provided for reference \cite{Tak8}.

\noindent \vspace{3pt}

In this manner, Claims (A) and (B) can be established in all cases. 

\vspace{10pt}

Now we explain how to prove Claim (C) following the strategy outlined
in Subsection \ref{subsec:(A) and (B) =00003D>(C)}. Similarly to
Claims (A) and (B), we check \eqref{0dimT} on the $p_{2}$-chart
and the locus $\{p_{1}=p_{2}=0\}|_{T}$ separately. On the $p_{2}$-chart
of $T_{\mA}$, the singularities are checked using the Jacobian criterion.
On $\{p_{1}=p_{2}=0\}|_{T}$, the verification has already been completed
since we have seen in the proofs of Claims (A) and (B) that this locus
consists of only finitely many points.

\noindent \vspace{3pt}

Hereafter in this subsection, the result of the computations for each
case follows.

\noindent \vspace{3pt}

\noindent \noindent%
\fbox{\textbf{No.$\,$360}}

\noindent \vspace{3pt}

\noindent \noindent\textbf{$\bullet$ Sections for $T$}:\begin{align*}& \begin{array}{|c|c|c|c|c|} \hline \text{weight} & 1 & 2 & 3 & 4 \\ \hline \text{section}  & s_{23}=0  & s_{13}=s_{22}=0  & s_{12}=t_3=0  & \begin{array}{c} q_3 = a_3 t_2 \\ s_{11} = a_{11} t_2 \end{array} \\ \hline \end{array} \\
& \begin{array}{|c|c|c|} \hline 5 & 6 & 8 \\ \hline q_2 = a_2 t_1 & q_1 = a_1 p_4 & u = a_0 p_2 + b_0 t_2^2 \\ \hline \end{array} 
\end{align*}

\noindent \noindent$\bullet$\textbf{ 6 parameters}: $a_{3},\dots,b_{0}\in\mC$

\noindent \noindent$\bullet$ \textbf{Embedding of $T$}:

\begin{equation}
T\subset\mP(t_{2},t_{1},p_{4},p_{1},r,p_{2},p_{3})=\mP(4,5,6,7^{2},8,9).\label{eq:360Tamb}
\end{equation}
\noindent$\bullet$ \textbf{Sing~$T$}: $T\cap\mP(t_{2},p_{2})$
consists of two points, say, $\mathsf{p}_{4},\mathsf{q}_{4}$. $T$
has $1/4(1,3)$-singularities at $\mathsf{p}_{4},\mathsf{q}_{4}$,
a $1/6(1,5)$-singularity at the $p_{4}$-point, and a $1/7(2,5)$-singularity
at the $p_{1}$-point.

\vspace{5pt}

\noindent %
\fbox{\textbf{No.$\,$393}}

\noindent \vspace{3pt}

\noindent \noindent\textbf{$\bullet$ Sections for $T$}: {\begin{align*} & \begin{array}{|c|c|c|c|c|} \hline \text{weight} & 1 & 2 & 3 & 4 \\ \hline \text{section} & s_{23} = 0 & \begin{array}{c} s_{13} =s_{22}= 0 \\  \end{array} & \begin{array}{c} s_{12} =t_3= 0 \\  \end{array} & \begin{array}{c} q_3 = a_3 t_2 \\ s_{11} = a_{11} t_2 \end{array} \\ \hline \end{array} \\ & \begin{array}{|c|c|c|} \hline 5 & 6 & 7 \\ \hline q_2 = a_2 p_4 + b_2 t_1 & q_1 = a_1 p_1 & r = a_0 p_2 \\ \hline \end{array} \end{align*} }

\noindent \noindent\textbf{$\bullet$ 6 parameters}:\textbf{ $a_{3},\dots,a_{0}\in\mC$} 

\noindent \noindent\textbf{$\bullet$ Embedding of $T$}:
\[
T\subset\mP(t_{2},p_{4},t_{1},p_{1},p_{2},p_{3},u)\simeq\mP(4,5^{2},6,7,8,9).
\]

\noindent \noindent$\bullet$ \textbf{Sing~$T$}: $T\cap\mP(p_{4},t_{1})$
consists of the $p_{4}$-point, and another point, say, $\mathsf{p}_{5}$.
$T$ has a $1/2(1,1)$-singularity at the point $\mathsf{p}_{2}:=T\cap\mP(t_{2},p_{1},p_{3})$
with $t_{2}p_{1}p_{3}\not=0$, a $1/5(1,4)$-singularity at the $p_{4}$-point,
a $1/5(2,3)$-singularity at $\mathsf{p}_{5}$, and\textbf{ }a $1/9(4,5)$-singularity
at the $u$-point (strictly speaking, $\mathsf{p}_{2}:=T\cap\mP(t_{2},p_{1},p_{3})$
should be written as $\{\mathsf{p}_{2}\}:=T\cap\mP(t_{2},p_{1},p_{3})$,
but we adopt this simplified notation here and in what follows).

\noindent \vspace{5pt}

\noindent %
\fbox{\textbf{No.$\,$569}}

\noindent \vspace{3pt}

\noindent \noindent$\bullet$\textbf{ Sections for $T$}: { \begin{align*} & \begin{array}{|c|c|c|c|c|} \hline \text{weight} & 1 & 2 & 3 & 4 \\ \hline \text{section} & s_{23} = 0 & \begin{array}{c} s_{13} =s_{22}= 0 \\  \end{array} & \begin{array}{c} p_4 = a_4 s_{12} \\ q_3 = a_3 s_{12} \end{array} & \begin{array}{c} q_2 = a_2 t_3 \\ s_{11} = a_{11} t_3 \end{array} \\ \hline \end{array} \\ & \begin{array}{|c|c|} \hline 5 & 6 \\ \hline \begin{array}{c} q_1 = a_1 p_1 + b_1 t_2 \end{array} & \begin{array}{c} r = a_0 p_2 + b_0 s_{12}^{2} \\ t_1 = c_1 p_2 + d_1 s_{12}^{2} \end{array} \\ \hline \end{array} \end{align*} } 

\noindent \noindent\textbf{$\bullet$ 10 parameters}:\textbf{ $a_{4},\dots,d_{1}\in\mC$} 

\noindent \noindent\textbf{$\bullet$ Embedding of $T$}:
\[
T\subset\mP(s_{12},t_{3},p_{1},t_{2},p_{2},p_{3},u)=\mP(3,4,5^{2},6,7,9).
\]

\noindent \noindent\textbf{$\bullet$ Sing~$T$}: $T\cap\mP(s_{12},p_{2},u)$
consists of the $u$-point and two points, say, $\sf{p}_{3},\sf{q}_{3}$.
$T$ has $1/3(1,2)$-singularities at $\sf{p}_{3},\sf{q}_{3}$, and
a $1/5(2,3)$-singularity at the point $\mathsf{p}_{5}:=T\cap\mP(p_{1},t_{2})$,
and a $1/9(4,5)$-singularity at the $u$-point. 

\noindent \vspace{5pt}

\noindent %
\fbox{\textbf{No.$\,$574}}

\noindent \vspace{3pt}

\noindent \noindent\textbf{$\bullet$ Sections for $T$}: {\begin{align*} & \begin{array}{|c|c|c|c|c|} \hline \text{weight} & 1 & 2 & 3 & 4 \\ \hline \text{section} & s_{23} = 0 & \begin{array}{c} s_{13} =s_{22}=t_3= 0  \\  \end{array} & \begin{array}{c} q_3 = a_3 t_2 \\ s_{12} = a_{12} t_2 \end{array} & \begin{array}{c} q_2 = a_2 t_1 \\ s_{11} = a_{11} t_1 \end{array} \\ \hline \end{array} \\ & \begin{array}{|c|c|} \hline 5 & 6 \\ \hline q_1 = a_1 p_1 + b_1 p_4 & r = a_0 p_2 + b_0 t_2^2 \\ \hline \end{array} \end{align*} }

\noindent \noindent\textbf{$\bullet$ 8 parameters}:\textbf{ $a_{3},\dots,b_{0}\in\mC$} 

\noindent \noindent\textbf{$\bullet$ Embedding of $T$}: 
\[
T\subset\mP(t_{2},t_{1},p_{1},p_{4},p_{2},p_{3},u)=\mP(3,4,5^{2},6,7^{2}).
\]

\noindent \noindent\textbf{$\bullet$ Sing~$T$}: $T\cap\mP(p_{1},p_{4})$
consists of the $p_{4}$-point and another point, say $\mathsf{p}_{5}$.
$T$ has a $1/3(1,2)$-singularity at the point $\mathsf{p}_{3}:=T\cap\mP(p_{2},t_{2})$,
a\textbf{ }$1/5(1,4)$-singularity at the $p_{4}$-point, a $1/5(2,3)$-singularity
at $\mathsf{p}_{5}$, and a $1/7(3,4)$-singularity at the $u$-point. 

\noindent \vspace{5pt}

\noindent %
\fbox{\textbf{No.$\,$642}}\,

\noindent \vspace{3pt}

\noindent \noindent\textbf{$\bullet$ Sections for $T$}: { \begin{align*} & \begin{array}{|c|c|c|c|c|} \hline \text{weight} & 1 & 2 & 3 & 4 \\ \hline \text{section} & s_{23} = 0 & \begin{array}{c} s_{13} =s_{22}= 0 \\ \end{array} & \begin{array}{c} p_4 = a_4 t_3 \\ s_{12} = a_{12} t_3 \end{array} & \begin{array}{c} \ s_{11} = a_{11} p_1 + b_{11} t_2 \\ q_3 = a_3 p_1 + b_3 t_2 \end{array} \\ \hline \end{array} \\ & \begin{array}{|c|c|} \hline 5 & 6 \\ \hline\begin{array}{c} q_2 = a_2 p_2\\ t_1=a_1 p_2\end{array} & \begin{array}{c} q_1 = b_1 p_3 + c_1 t_3^2\\ \end{array} \\ \hline \end{array} \end{align*} }

\noindent \noindent\textbf{$\bullet$ 10 parameters}:\textbf{ $a_{4},\dots,c_{1}\in\mC$} 

\noindent \noindent\textbf{$\bullet$ Embedding of $T$}: 
\[
T\subset\mP(t_{3},p_{1},t_{2},p_{2},p_{3},r,u)=\mP(3,4^{2},5,6,7,11).
\]
\noindent\textbf{$\bullet$ Sing~$T$}: $T$ has a $1/2(1,1)$-singularity
at the point 
\[
\mathsf{p}_{2}:=T\cap\mP(p_{1},t_{2},p_{3})\cap\{(a_{3}+b_{1})p_{1}+b_{3}t_{2}=0\},
\]

\noindent a $1/3(1,2)$-singularity at the point ${\sf{p}_{3}}:=T\cap\mP(p_{3},t_{3}),$
a $1/4(1,3)$-singularity at the point ${\sf{p}_{4}}:=T\cap\mP(p_{1},t_{2})$\textbf{,
}and a $1/11(4,7)$-singularity at the $u$-point.

\noindent \vspace{5pt}

\noindent %
\fbox{\textbf{No.$\,$644}}

\noindent \vspace{3pt}

\noindent \noindent\textbf{$\bullet$ Sections for $T$}: { \begin{align*} & \begin{array}{|c|c|c|c|c|} \hline \text{weight} & 1 & 2 & 3 & 4 \\ \hline \text{section} & s_{23} = 0 & \begin{array}{c} s_{13} =s_{22}=t_3 = 0 \\  \end{array} & s_{12} = a_{12} t_2 & \begin{array}{c} p_4 = a_4 p_1 + b_4 t_1 \\ q_3 = a_3 p_1 + b_3 t_1 \\ \ s_{11} = a_{11} p_1 + b_{11} t_1 \end{array} \\ \hline \end{array} \\ & \begin{array}{|c|c|} \hline 5 & 6 \\ \hline q_2 = a_2 p_2 & q_1 = a_1 p_3 + b_1 t_2^2 \\ \hline \end{array} \end{align*} } 

\noindent \noindent\textbf{$\bullet$ 10 parameters}:\textbf{ $a_{12},\dots,b_{1}\in\mC$} 

\noindent \noindent\textbf{$\bullet$ Embedding of $T$}:
\[
T\subset\mP(t_{2},p_{1},t_{1},p_{2},p_{3},r,u)=\mP(3,4^{2},5,6,7,10)
\]

\noindent \noindent\textbf{$\bullet$ Sing~$T$}: $T\cap\mP(p_{1},t_{1})$
consists of two points, say, ${\sf p}_{4}$, ${\sf q}_{4}$. $T$
has a\textbf{ }$1/2(1,1)$-singularity at the point 
\[
{\sf p}_{2}:=T\cap\mP(p_{1},t_{1},p_{3})\cap\{(a_{1}+a_{3})p_{1}+b_{3}t_{1}=0\},
\]
$1/4(1,3)$-singularities at $\mathsf{p}_{4}$ and $\mathsf{q}_{4}$,
and a\textbf{ }$1/10(3,7)$-singularity at the $u$-point.

\noindent \vspace{5pt}

\noindent %
\fbox{\textbf{No.$\,$1091}}

\noindent \vspace{3pt}

\noindent \noindent\textbf{$\bullet$ Sections for $T$}: { \begin{align*} & \begin{array}{|c|c|c|c|c|} \hline \text{weight} & 1 & 2 & 3 & 4 \\ \hline \text{section} & s_{23} = 0 & s_{22} = a_{22} s_{13} & \begin{array}{c} s_{12} =t_3= 0 \\  \end{array} & \begin{array}{c} q_3 = a_3 s_{13}^{2} \\ s_{11} = a_{11} s_{13}^{2} \\ t_2 = a_2 s_{13}^{2} \end{array} \\ \hline \end{array} \\ & \begin{array}{|c|c|} \hline 5 & 6 \\ \hline \begin{array}{c} p_4 = a_4 t_1 \\ q_2 = b_2 t_1 \end{array} & q_1 = a_1 p_1 + b_1 s_{13}^{3} \\ \hline \end{array} \end{align*} } 

\noindent \noindent\textbf{$\bullet$ 8 parameters}:\textbf{ $a_{22},\dots,b_{1}\in\mC$} 

\noindent \noindent\textbf{$\bullet$ Embedding of $T$}: 
\[
T\subset\mP(s_{13},t_{1},p_{1},r,p_{2},p_{3},u)=\mP(2,5,6,7^{2},8,9).
\]

\noindent \noindent\textbf{$\bullet$ Sing~$T$}: $T\cap\mP(s_{13},p_{1},p_{3})$
consists of two points, say, $\sf{p}_{2},\sf{q}_{2}$ with $s_{13}\not=0$.
$T$ has\textbf{ }$1/2(1,1)$-singularities at ${\sf p}_{2}$ and
${\sf q}_{2}$, a $1/7(1,6)$-singularity at the $p_{2}$-point, and
a\textbf{ }$1/9(2,7)$-singularity at the $u$-point.

\noindent \vspace{5pt}

\noindent %
\fbox{\textbf{No.$\,$1181}}

\noindent \vspace{3pt}

\noindent \noindent\textbf{$\bullet$ Sections for $T$} { \begin{align*} & \begin{array}{|c|c|c|c|c|} \hline \text{weight} & 1 & 2 & 3 & 4 \\ \hline \text{section} & s_{23} = 0 & \begin{array}{c} p_4 = a_4 s_{13} \\ s_{22} = a_{22} s_{13} \end{array} & \begin{array}{c} t_3 = a_3 p_1 \\ s_{12} = a_{12} p_1 \end{array} & \begin{array}{c} \ s_{11} = a_{11} p_2 + b_{11} s_{13}^2 \\ q_3 = b_3 p_2 + c_3 s_{13}^2 \\ t_2 = a_2 p_2 + b_2 s_{13}^2 \end{array} \\ \hline \end{array} \\ & \begin{array}{|c|c|} \hline 5 & 6 \\ \hline q_2 = c_2 p_3 + d_2 t_1 + e_2 s_{13} p_1 & q_1 = a_1 p_2 s_{13} + b_1 s_{13}^3 + c_1 p_1^2 \\ \hline \end{array} \end{align*} }

\noindent \noindent\textbf{$\bullet$ 16 parameters $a_{4},\dots,c_{1}\in\mC$} 

\noindent \noindent\textbf{$\bullet$ Embedding of $T$}: 
\[
T\subset\mP(s_{13},p_{1},p_{2},p_{3},t_{1},r,u)=\mP(2,3,4,5^{2},7,12).
\]

\noindent \noindent\textbf{$\bullet$ Sing~$T$}: $T\cap\mP(p_{2},u)$
consists of the $u$-point and another point, say, ${\sf p}_{4}$.
$T$ has a $1/2(1,1)$-singularity at the point $\mathsf{p}_{2}:=T\cap\mP(s_{13},p_{2},u)\cap\{s_{13}\not=0\}$,
a\textbf{ }$1/4(1,3)$-singularity at ${\sf p}_{4}$, and a $1/12(5,7)$-singularity
at the $u$-point. 

\noindent \vspace{5pt}

\noindent %
\fbox{\textbf{No.$\,$1185}}

\noindent \vspace{3pt}

\noindent \noindent\textbf{$\bullet$ Sections for $T$}: { \begin{align*} & \begin{array}{|c|c|c|c|c|} \hline \text{weight} & 1 & 2 & 3 & 4 \\ \hline \text{section} & s_{23} = 0 & \begin{array}{c} p_4 = a_4 s_{13} \\ q_3 = a_3 s_{13} \\ s_{22} = a_{22} s_{13} \end{array} & q_2 = a_2 s_{12} & \begin{array}{c} q_1 = a_1 p_1 + b_1 s_{13}^{2} \\ \ s_{11} = a_{11} p_1 + b_{11} s_{13}^{2} \\ t_3 = c_3 p_1 + d_3 s_{13}^{2} \end{array} \\ \hline \end{array} \\ & \begin{array}{|c|c|} \hline 5 & 6 \\ \hline r = a_0 p_2 + b_0 t_2 + c_0 s_{13} s_{12} & t_1 = c_1 p_3 + d_1 s_{12}^2 + e_1 s_{13}^3 + f_1 p_1 s_{13} \\ \hline \end{array} \end{align*} } 

\noindent \noindent\textbf{$\bullet$ 17 parameters}:\textbf{ $a_{4},\dots,f_{1}\in\mC$} 

\noindent \noindent\textbf{$\bullet$ Embedding of $T$}:
\begin{equation}
T\subset\mP(s_{13},s_{12},p_{1},p_{2},t_{2},p_{3},u)=\mP(2,3,4,5^{2},6,8).\label{eq:1185X}
\end{equation}
\noindent\textbf{$\bullet$ Sing~$T$}: $T\cap\mP(s_{13},p_{1},p_{3},u)\cap\{s_{13}\not=0\}$
consists of three points, say, $\mathsf{p}_{2}$, $\mathsf{q}_{2}$,
$\mathsf{r}_{2}$. $T$ has $1/2(1,1)$-singularities at $\mathsf{p}_{2}$,
$\mathsf{q}_{2}$, $\mathsf{r}_{2}$, a\textbf{ }$1/5(1,4)$-singularity
at the point ${\sf{p}_{5}}:=T\cap\mP(p_{2},t_{2}),$ and a\textbf{
}$1/8(3,5)$-singularity at the $u$-point. 

\noindent \noindent$\bullet$ \textbf{Comments: }We describe how
to determine the singularities of $T$ along the locus $\{p_{1}=p_{2}=0\}|_{T}$.
The detailed computations are provided in the accompanying \textit{Mathematica}
code \cite{Tak8}. We see that \textbf{$\{p_{1}=p_{2}=0\}_{|T}$ }consists
of the $u$-point and the locus $\Delta:=T\cap\mP(s_{13},s_{12},t_{2},u)\cap\{s_{13}\not=0\}$,
and $\Delta$ consists of two points. Applying LPC at the $u$-point,
we see that $T$ has a $1/8(3,5)$-singularity there. We show that
$T$ is smooth along $\Delta$. We observe that none of the coordinates
$s_{13},s_{12},t_{2},u$ of the points in $\Delta$ are zero. In particular,
$\Delta$ is contained in the $u$-chart of $T$. We can check that
the $u$-chart of $W_{\mA}:=\SigmaA\cap\{s_{23}=0\}$ is the affine
$12$-space. Since the weights of $s_{13},s_{12},t_{2},u$ are $2,3,5,8$,
respectively, the points in $\Delta$ are not contained in $\Bs|\sO_{\mP_{\Sigma}}(i)|$
with $i=2,3,4,5,6$. Hence $T_{\mA}^{o}$ is smooth at the point in
the inverse image of $\Delta$ by $T_{\mA}^{o}\to T$ since the $u$-chart
of $W_{\mA}$ is smooth and $T_{\mA}=W_{\mA}\cap(2)^{3}\cap(3)\cap(4)^{3}\cap(5)\cap(6)$.
Since $w(u)=8$, the cyclic group $\mZ_{8}$ acts on ${\mA}_{T}^{u}$.
Since none of the coordinates $s_{13},s_{12},t_{2},u$ of the point
in $\Delta$ are zero, the $\mZ_{8}$-action on $\mA_{T}^{u}$ is
free along the inverse image of $\Delta$. Therefore $T$ is smooth
along $\Delta$.

\noindent \vspace{5pt}

\noindent %
\fbox{\textbf{No.$\,$1186}}

\noindent \vspace{3pt}

\noindent \noindent\textbf{$\bullet$ Sections for $T$}: {\begin{align*} & \begin{array}{|c|c|c|c|c|} \hline \text{weight} & 1 & 2 & 3 & 4 \\ \hline \text{section} & s_{23} = 0 & \begin{array}{c} q_3 = a_3 s_{13} \\ s_{22} = a_{22} s_{13} \end{array} & \begin{array}{c} p_4 = a_4 t_3 \\ q_2 = a_2 t_3 \\ s_{12} = a_{12} t_3 \end{array} & \begin{array}{c} q_1 = a_1 p_1 + b_1 s_{13}^{2} \\       \ \, s_{11} = a_{11} p_1 + b_{11} s_{13}^{2} \\ t_2 = b_2 p_1 + c_2 s_{13}^{2} \end{array} \\ \hline \end{array} \\ & \begin{array}{|c|} \hline 5 \\ \hline r = a_0 p_2 + b_0 t_1 + c_0 s_{13} t_3 \\ \hline \end{array} \end{align*} } 

\noindent \noindent\textbf{$\bullet$ 14 parameters}:\textbf{ $a_{3},\dots,c_{0}\in\mC$} 

\noindent \noindent\textbf{$\bullet$ Embedding of $T$}:
\[
T\subset\mP(s_{13},t_{3},p_{1},p_{2},t_{1},p_{3},u)=\mP(2,3,4,5^{2},6,7).
\]
\noindent\textbf{$\bullet$ Sing~$T$}: $T\cap\mP(s_{13},p_{1},p_{3})$
consists of two points with $s_{13}\not=0$, say, $\mathsf{p}_{2}$,
$\mathsf{q}_{2}$. $T$ has $1/2(1,1)$-singularities at $\mathsf{p}_{2}$,
$\mathsf{q}_{2}$, a\textbf{ }$1/3(1,2)$-singularity at the point
$\mathsf{p}_{3}:=T\cap\mP(t_{3},p_{3})$, a\textbf{ }$1/5(1,4)$-singularity
at the point ${\sf{p}_{5}}:=T\cap\mP(p_{2},t_{1}),$ and a\textbf{
}$1/7(2,5)$-singularity at the $u$-point.

\noindent \vspace{5pt}

\noindent %
\fbox{\textbf{No.$\,$1218}}

\noindent \vspace{3pt}

\noindent \noindent\textbf{$\bullet$ Sections for $T$}: { \begin{align*} & \renewcommand{\arraystretch}{1.4} \begin{array}{|c|c|c|c|c|} \hline \text{weight} & 1 & 2 & 3 & 4 \\ \hline \text{section} & s_{23} = t_{3} = 0 & \begin{array}{c} q_3 = a_3 t_2 \\ s_{13} = a_{13} t_2 \\ s_{22} = a_{22} t_2 \end{array} & \begin{array}{c} q_2 = a_2 t_1 \\ s_{12} = a_{12} t_1 \end{array} & \begin{array}{c} q_1 = a_1 p_1 + b_1 t_2^2 \\ \ \, s_{11} = a_{11} p_1 + b_{11} t_2^2 \end{array} \\ \hline \end{array} \\[1ex] & \begin{array}{|c|} \hline 5 \\ \hline p_4 = a_4 p_2 + b_4 r + c_4 u + d_4 t_2 t_1 \\ \hline \end{array} \end{align*} } 

\noindent \noindent\textbf{$\bullet$ 13 parameters}:\textbf{ $a_{3},\dots,d_{4}\in\mC$} 

\noindent \noindent\textbf{$\bullet$ Embedding of $T$}:
\[
T\subset\mP(t_{2},t_{1},p_{1},p_{2},r,u,p_{3})=\mP(2,3,4,5^{3},6).
\]

\noindent \noindent\textbf{$\bullet$ Sing~$T$}: $T\cap\mP(t_{2},p_{1},p_{3})$
consists of two points with $t_{2}\not=0$, say, $\mathsf{p}_{2}$,
$\mathsf{q}_{2}$. $T\cap\mP(p_{2},r,u)$ consists of two points,
say, $\mathsf{p}_{5}$, $\mathsf{q}_{5}$. $T$ has $1/2(1,1)$-singularities
at $\mathsf{p}_{2}$, $\mathsf{q}_{2}$, $1/5(2,3)$-singularities
at $\mathsf{p}_{5}$, $\mathsf{q}_{5}$, and a $1/5(1,4)$-singularity
at the $p_{2}$-point. 

\noindent \vspace{5pt}

\noindent %
\fbox{\textbf{No.$\,$1253}}

\noindent \vspace{3pt}

\noindent \noindent\textbf{$\bullet$ Sections for $T$}: { \begin{align*} &\renewcommand{\arraystretch}{1.4} \begin{array}{|c|c|c|} \hline \text{weight} & 1 & 2 \\ \hline \text{section} & s_{13} = s_{23} = 0 & \begin{array}{c} s_{11} = a_{11} t_{3} \\ s_{12} = a_{12} t_{3} \\ s_{22} = a_{22} t_{3} \end{array} \\ \hline \end{array} \\[1ex] &\renewcommand{\arraystretch}{1.4} \begin{array}{|c|c|} \hline 3 & 4 \\ \hline \begin{array}{c} p_{4} = a_{4} t_{2} \\ q_{3} = a_{3} t_{2} \\ t_{1} = a_{1} t_{2} \end{array} & \begin{array}{c} q_{1} = b_{1} p_{1} + c_{1} p_{2} + d_{1} t_{3}^{2} \\ q_{2} = a_{2} p_{1} + b_{2} p_{2} + c_{2} t_{3}^{2} \end{array} \\ \hline \end{array} \end{align*} } 

\noindent \noindent\textbf{$\bullet$ 12 parameters}:\textbf{ $a_{11},\dots,c_{2}\in\mC$} 

\noindent \noindent\textbf{$\bullet$ Embedding of $T$}:
\[
T\subset\mP(t_{3},t_{2},p_{1},p_{2},p_{3},r,u)=\mP(2,3,4^{2},5^{2},7).
\]

\noindent \noindent\textbf{$\bullet$ Sing~$T$}: $T\cap\mP(t_{3},p_{1},p_{2})\cap\{t_{3}\not=0\}$
consists of two points, say, $\mathsf{p}_{2}$, $\mathsf{q}_{2}$.
$T\cap\mP(p_{1},p_{2})$ consists of two points, say, $\mathsf{p}_{4}$,
$\mathsf{q}_{4}$. $T$ has $1/2(1,1)$-singularities at $\mathsf{p}_{2}$,
$\mathsf{q}_{2}$, $1/4(1,3)$-singularities at $\mathsf{p}_{4}$,
$\mathsf{q}_{4}$, and a\textbf{ }$1/7(2,5)$-singularity at the $u$-point.

\noindent \vspace{5pt}

\noindent %
\fbox{\textbf{No.$\,$1413}}

\noindent \vspace{3pt}

\noindent \noindent\textbf{$\bullet$ Sections for $T$}: { \begin{align*} &\renewcommand{\arraystretch}{1.4} \begin{array}{|c|c|c|} \hline \text{weight} & 1 & 2 \\ \hline \text{section} & s_{13} = s_{23} = 0 & \begin{array}{c} q_{3} = a_{3} t_{3} \\ s_{11} = a_{11} t_{3} \\ s_{12} = a_{12} t_{3} \\ s_{22} = a_{22} t_{3} \end{array} \\ \hline \end{array} \\[1ex] &\renewcommand{\arraystretch}{1.4} \begin{array}{|c|c|} \hline 3 & 4 \\ \hline \begin{array}{c} p_{4} = a_{4} t_{1} + b_{4} t_{2} \\ q_{1} = a_{1} t_{1} + b_{1} t_{2} \\ q_{2} = a_{2} t_{1} + b_{2} t_{2} \end{array} & r = a_{0} p_{1} + b_{0} p_{2} + c_{0} t_{3}^{2} \\ \hline \end{array} \end{align*} } 

\noindent \noindent\textbf{$\bullet$ 13 parameters $a_{3},\dots,c_{0}\in\mC$} 

\noindent \noindent\textbf{$\bullet$ Embedding of $T$}:
\[
T\subset\mP(t_{3},t_{1},t_{2},p_{1},p_{2},p_{3},u)=\mP(2,3^{2},4^{2},5^{2}).
\]
\noindent$\bullet$\textbf{ Sing~$T$}: $T\cap\mP(t_{3},p_{1},p_{2})\cap\{t_{3}\not=0\}$
consists of two points, say, $\mathsf{p}_{2}$, $\mathsf{q}_{2}$.
$T\cap\mP(t_{1},t_{2})$ consists of two points, say, $\mathsf{p}_{3}$,
$\mathsf{q}_{3}$. $T$ has $1/2(1,1)$-singularities at $\mathsf{p}_{2}$,
$\mathsf{q}_{2}$, $1/3(1,2)$-singularities at $\mathsf{p}_{3}$,
$\mathsf{q}_{3}$, a\textbf{ }$1/4(1,3)$-singularity at the point
$\mathsf{p}_{4}:=T\cap\mP(p_{1},p_{2})$, and a\textbf{ }$1/5(2,3)$-singularity
at the $u$-point.

\noindent \vspace{5pt}

\noindent %
\fbox{\textbf{No.$\,$2422}}

\noindent \vspace{3pt}

\noindent \noindent\textbf{$\bullet$ Sections for $T$}: { \begin{align*} & \renewcommand{\arraystretch}{1.4} \begin{array}{|c|c|c|c|} \hline \text{weight} & 1 & 2 & 3 \\ \hline \text{section} & t_{3} = s_{23} = 0 & \begin{array}{c} q_{3} = a_{3} p_{1} + b_{3} t_{2} \\ \ \, s_{13} = a_{13} p_{1} + b_{13} t_{2} \\ \ \, s_{22} = a_{22} p_{1} + b_{22} t_{2} \end{array} & \begin{array}{c} p_{4} = a_{4} p_{2} + b_{4} t_{1} \\ q_{2} = a_{2} p_{2} + b_{2} t_{1} \\ \ \, s_{12} = a_{12} p_{2} + b_{12} t_{1} \end{array} \\ \hline \end{array} \\[1ex] & \begin{array}{|c|} \hline  4 \\ \hline  \begin{array}{c} q_{1} = a_{1} p_{3} + b_{1} p_{1}^{2} + c_{1} p_{1} t_{2} + d_{1} t_{2}^{2} \\ \quad \, s_{11} = a_{11} p_{3} + b_{11} p_{1}^{2} + c_{11} p_{1} t_{2} + d_{11} t_{2}^{2} \end{array} \\ \hline \end{array} \end{align*} \noindent\textbf{$\bullet$
20 parameters}:\textbf{ $a_{3},\dots,d_{11}\in\mC$} 

\noindent \noindent\textbf{$\bullet$ Embedding of $T$}:
\[
T\subset\mP(p_{1},t_{2},p_{2},t_{1},p_{3},r,u)=\mP(2^{2},3^{2},4,5,7).
\]

\noindent \noindent$\bullet$\textbf{ Sing~$T$}: $T\cap\mP(p_{1},t_{2},p_{3})$
consists of four points with $p_{1}\not=0$, say, $\mathsf{p}_{2},\mathsf{q}_{2},\mathsf{r}_{2},\mathsf{s}_{2}$.
$T$ has $1/2(1,1)$-singularities at $\mathsf{p}_{2},\mathsf{q}_{2},\mathsf{r}_{2},\mathsf{s}_{2}$,
a\textbf{ }$1/3(1,2)$-singularity at the point ${\sf{p}_{3}}:=T\cap\mP(p_{2},t_{1})$,
and a\textbf{ }$1/7(2,5)$-singularity at the $u$-point.

\subsection{Case $h^{0}(\sO_{\mP_{X}}(1))=2$\label{subsec:Caseh0=00003D2}}

\noindent We outline the proofs of Claims (A) and (B). Let $C$ be
as in Subsection \ref{subsec:Presentation-of-the eq XT}. To show
Claim (A) for $T$, it suffices to show that $C_{\mA}^{0}$ is smooth.
Indeed, this reduction can be shown by noting $h^{0}(\sO_{\mP_{T}}(1))=1$
and repeating for $T$ and $C$ the argument for $X$ and $T$ in
Subsection \ref{subsec:Reduction-of-Claims}. The smoothness of $C_{\mA}^{0}$
is proved separately on the $p_{1}$-chart, the $p_{2}$-chart, and
the locus $\{p_{1}=p_{2}=0\}_{|C}$ for No.$\,$5866, 6860, and 6865
(resp.$\,$on one chosen coordinate chart and its complement for No.$\,$4850,
4938, 5202, and 5859). The smoothness of $C_{\mA}^{0}$ can be proved
on the $p_{1}$- and $p_{2}$-charts (resp.$\,$on one chosen coordinate
chart) by the Jacobian criterion. It turns out that the locus $\{p_{1}=p_{2}=0\}_{|C}$
(resp.$\,$the complement of the coordinate chart) consists of only
the $u$-point. By LPC, we can determine the singularity of $T$ at
the $u$-point (cf.$\,$Subsection \ref{subsec:LPC-for C} and Example
\ref{exa:LPC at u-pt}). Thus Claim (A) for $T$ follows, and Claim
(B) for $T$ is also proved at the $u$-point. It remains to show
Claim (B) for $T$ outside the $u$-point. Similarly to the proof
of Claim \ref{cla:RedAB}, we see that $\Sing T$ is contained in
$C$. Thus $T$ has singularities along the quotient of the nonfree
locus of the $\mC^{*}$-action on $C_{\mA}^{o}$. Consequently, for
the classes except No.$\,$4938, we prove that $T$ has only $1/2(1,1)$-singularities
on the $p_{1}$- and $p_{2}$-chart (resp.$\,$the chosen coordinate
chart) following the argument in Subsection \ref{subsec:1/2-singularities};
for No.$\,$4938, we prove that $T$ has a 1/3(1,2)-singularity on
the $p_{1}$-chart by LPC.

Now we outline the proof of Claim (C). Similarly to the argument in
Subsection \ref{subsec:(A) and (B) =00003D>(C)}, we see that it suffices
to show $\dim\Sing(C\cap\{p_{1}=0\})\leq0$. We can check that $C\cap\{p_{1}=0\}$
itself consists of a finite number of points except in the case of
No.$\,$5859. Thus the claim holds in all the other cases. For No.$\,$5859,
it holds that $C\subset\{p_{1}=0\}$, hence the previous argument
does not apply. In this case, we change the definition of $\mathsf{b}$.
By the ${\rm GL}_{3}$-action on $\Sigma_{\mA}^{14}$ (\cite[Prop.\,2.11]{Tak9}),
we may assume $\mathsf{b}=p_{2}$. Then $C\cap\{p_{2}=0\}$ itself
consists of a finite number of points, and the claim follows.

\noindent \vspace{5pt}

\noindent %
\fbox{\textbf{No.$\,$4850}}

\noindent \vspace{3pt}

\noindent \noindent\textbf{$\bullet$ Sections for $C$} { 
\begin{center}
{\small{}\begin{align*} & \renewcommand{\arraystretch}{1.4} \begin{array}{|c|c|c|c|c|} \hline \text{weight} & 1 & 2 & 3 & 4 \\ \hline \text{section} & p_{4} = s_{23} = 0 & s_{13} = s_{22} = 0 & s_{12} = 0 & \begin{array}{c} s_{11} = a_{11} p_{1} \\ q_{3} = a_{3} p_{1} \end{array} \\ \hline \end{array} \\[1ex] & \begin{array}{|c|c|c|} \hline  5 & 6 & 7 \\ \hline  \begin{array}{c} q_{2} = a_{2} p_{2} \\ t_{3} = b_{3} p_{2} \end{array} & q_{1} = a_{1} p_{3} + b_{1} t_{2} & r = a_{0} t_{1} \\ \hline \end{array} \end{align*}  }{\small\par}
\par\end{center}

} 

\noindent\textbf{$\bullet$ 7 parameters}:\textbf{ $a_{11},\dots,a_{0}\in\mC$} 

\noindent\textbf{$\bullet$ Embedding of $C$}:
\[
C\subset\mP(p_{1},p_{2},p_{3},t_{2},t_{1},u)=\mP(4,5,6^{2},7,13).
\]

\noindent \noindent$\bullet$\textbf{ Sing~$T$}: $T$ has a $1/2(1,1)$-singularity
at the point $\mathsf{p}_{2}:=T\cap\mP(p_{1},p_{3},t_{2})$ with $p_{1}p_{3}t_{2}\not=0$,
and a $1/13(6,7)$-singularity at the $u$-point. 

\noindent \vspace{5pt}

\noindent %
\fbox{\textbf{No.$\,$4938}}

\noindent \vspace{3pt}

\noindent \noindent\textbf{$\bullet$ Sections for $C$}: { \begin{align*} & \renewcommand{\arraystretch}{1.4} \begin{array}{|c|c|c|c|c|} \hline \text{weight} & 1 & 2 & 3 & 4 \\ \hline \text{section} & \begin{array}{c} p_{4} =s_{23}=0 \\  \end{array} & \begin{array}{c} s_{13} =s_{22}=0 \\  \end{array} & \begin{array}{c} s_{12} = a_{12} p_{1} \\ q_{3} = a_{3} p_{1} \end{array} & \begin{array}{c} q_{2} = a_{2} p_{2} \\ s_{11} = a_{11} p_{2} \\ t_{3} = b_{3} p_{2} \end{array} \\ \hline \end{array} \\[1ex] & \renewcommand{\arraystretch}{1.4} \begin{array}{|c|c|} \hline  5 & 6 \\ \hline  q_{1} = a_{1} p_{3} + b_{1} t_{2} & r = a_{0} t_{1} + b_{0} p_{1}^2 \\ \hline \end{array} \end{align*} } 

\noindent \noindent\textbf{$\bullet$ 9 parameters}:\textbf{ $a_{12},\dots,b_{0}\in\mC$} 

\noindent \noindent\textbf{$\bullet$ Embedding of $C$}: 
\[
C\subset\mP(p_{1},p_{2},p_{3},t_{2},t_{1},u)=\mP(3,4,5^{2},6,11).
\]

\noindent \noindent$\bullet$\textbf{ Sing~$T$}: $T$ has a $1/3(1,2)$-singularity
at the point $\mathsf{p}_{3}:=T\cap\mP(p_{1},t_{1})$, and a $1/11(5,6)$-singularity
at the $u$-point. 

\noindent \vspace{5pt}

\noindent %
\fbox{\textbf{No.$\,$5202}}

\noindent \vspace{3pt}

\noindent \noindent\textbf{$\bullet$ Sections for $C$}: { 
\begin{center}
{\small{}\begin{align*} &\renewcommand{\arraystretch}{1.4} \begin{array}{|c|c|c|c|} \hline \text{weight} & 1 & 2 & 3 \\ \hline \text{section} & p_{4} = s_{23} = 0 & \begin{array}{c} q_{3} = a_{3} p_{1} \\ s_{22} = a_{22} p_{1} \\ s_{13} = a_{13} p_{1} \end{array} & \begin{array}{c} q_{2} = a_{2} p_{2} \\ s_{12} = a_{12} p_{2} \\ t_{3} = a_{0} p_{2} \end{array} \\ \hline \end{array} \\[1ex] &\renewcommand{\arraystretch}{1.4} \begin{array}{|c|c|} \hline 4 & 5 \\ \hline \begin{array}{c} q_{1} = a_{1} p_{3} + b_{1} t_{2} + c_{1} p_{1}^{2} \\ \ \ \ s_{11} = a_{11} p_{3} + b_{11} t_{2} + c_{11} p_{1}^{2} \end{array} & r = a_{0} t_{1} + b_{0} p_{1} p_{2} \\ \hline \end{array} \end{align*}  }{\small\par}
\par\end{center}

} 

\noindent\textbf{$\bullet$ 14 parameters}:\textbf{ $a_{3},\dots,b_{0}\in\mC$} 

\noindent\textbf{$\bullet$ Embedding of $C$}: 
\[
C\subset\mP(p_{1},p_{2},p_{3},t_{2},t_{1},u)=\mP(2,3,4^{2},5,9).
\]

\noindent \noindent$\bullet$\textbf{ Sing~$T$}: $T\cap\mP(p_{1},p_{3},t_{2})$
consists of two points with $p_{1}\not=0$, say, $\mathsf{p_{2}},\mathsf{q}_{2}$.
$T$ has $1/2(1,1)$-singularities at $\mathsf{p_{2}},\mathsf{q}_{2}$,
and a $1/9(4,5)$-singularity at the $u$-point. 

\noindent \vspace{5pt}

\noindent %
\fbox{\textbf{No.$\,$5859}}

\noindent \vspace{3pt}

\noindent \noindent\textbf{$\bullet$ Sections for $C$}: { 
\begin{center}
{\small{}\begin{align*} &\renewcommand{\arraystretch}{1.4} \begin{array}{|c|c|c|} \hline \text{weight} & 1 & 2 \\ \hline \text{section} & t_{3} = p_{1} = s_{23} = 0 & \begin{array}{c} p_{4} = a_{4} p_{2} + b_{4} t_{2} \\ q_{3} = a_{3} p_{2} + b_{3} t_{2} \\ \ \, s_{13} = a_{13} p_{2} + b_{13} t_{2} \\ \ \, s_{22} = a_{22} p_{2} + b_{22} t_{2} \end{array} \\ \hline \end{array} \\[1ex] &\renewcommand{\arraystretch}{1.4} \begin{array}{|c|c|} \hline 3 & 4 \\ \hline \begin{array}{c} q_{2} = a_{2} p_{3} + b_{2} t_{1} \\ \ \, s_{12} = a_{12} p_{3} + b_{12} t_{1} \end{array} & \begin{array}{c} q_{1} = a_{1} p_{2}^{2} + b_{1} p_{2} t_{2} + c_{1} t_{2}^{2} \\ \ \ \ s_{11} = a_{11} p_{2}^{2} + b_{11} p_{2} t_{2} + c_{11} t_{2}^{2} \end{array} \\ \hline \end{array} \end{align*}  }{\small\par}
\par\end{center}

} 

\noindent\textbf{$\bullet$ 18 parameters}:\textbf{ $a_{3},\dots,b_{0}\in\mC$} 

\noindent\textbf{$\bullet$ Embedding of $C$}: 
\[
C\subset\mP(p_{2},t_{2},p_{3},t_{1},r,u)=\mP(2^{2},3^{2},5,8).
\]

\noindent \noindent$\bullet$\textbf{ Sing~$T$}: $\mP(p_{2},t_{2},u)\cap T$
consists of the $u$-point, and other two points with $p_{2}\not=0$,
say, $\mathsf{p}_{2},\mathsf{q}_{2}$. $T$ has $1/2(1,1)$-singularities
at $\mathsf{p}_{2},\mathsf{q}_{2}$, and a $1/8(3,5)$-singularity
at the $u$-point.

\noindent \vspace{5pt}

\noindent %
\fbox{\textbf{No.$\,$5866}}

\noindent \vspace{3pt}

\noindent \noindent\textbf{$\bullet$ Sections for $C$}: { 
\begin{center}
\begin{align*} &\renewcommand{\arraystretch}{1.4} \begin{array}{|c|c|c|} \hline \text{weight} & 1 & 2 \\ \hline \text{section} & p_{4} = s_{13} = s_{23} = 0 & \begin{array}{c} q_{3} = a_{3} p_{1} + b_{3} p_{2} \\ \ s_{11} = a_{11} p_{1} + b_{11} p_{2} \\ \ s_{12} = a_{12} p_{1} + b_{12} p_{2} \\ \ s_{22} = a_{22} p_{1} + b_{22} p_{2} \\ t_{3} = a_{0} p_{1} + b_{0} p_{2} \end{array} \\ \hline \end{array} \\[1ex] &\renewcommand{\arraystretch}{1.4} \begin{array}{|c|} \hline 3 \\ \hline \begin{array}{c} q_{1} = a_{1} p_{3} + b_{1} t_{2} \\ q_{2} = a_{2} p_{3} + b_{2} t_{2} \\ t_{1} = c_{1} p_{3} + d_{1} t_{2} \end{array} \\ \hline \end{array} \end{align*} {\small{} }{\small\par}
\par\end{center}

} 

\noindent\textbf{$\bullet$ 16 parameters}:\textbf{ $a_{3},\dots,d_{1}\in\mC$} 

\noindent\textbf{$\bullet$ Embedding of $C$}: 
\[
C\subset\mP(p_{1},p_{2},p_{3},t_{2},r,u)=\mP(2^{2},3^{2},4,7).
\]

\noindent \noindent$\bullet$\textbf{ Sing~$T$}: $\mP(p_{1},p_{2},r)\cap T$
consists of three points with $p_{1}\not=0$, say, $\mathsf{p}_{2},\mathsf{q}_{2},\mathsf{r}_{2}$.
$T$ has $1/2(1,1)$-singularities at $\mathsf{p}_{2},\mathsf{q}_{2}$,
$\mathsf{r}_{2}$ and a $1/7(3,4)$-singularity at the $u$-point.

\noindent \vspace{5pt}

\noindent %
\fbox{\textbf{No.$\,$6860}}

\noindent \vspace{3pt}

\noindent \noindent\textbf{$\bullet$ Sections for $C$}: { \begin{align*} &\renewcommand{\arraystretch}{1.4} \begin{array}{|c|c|c|} \hline \text{weight} & 1 & 2 \\ \hline \text{section} & p_{4} = q_{3} = s_{13} = s_{23} = 0 & \begin{array}{c} \ \ \, s_{11} = a_{11} p_{1} + b_{11} p_{2} + c_{11} t_{3} \\ \ \ \, s_{12} = a_{12} p_{1} + b_{12} p_{2} + c_{12} t_{3} \\ \ \ \, s_{22} = a_{22} p_{1} + b_{22} p_{2} + c_{22} t_{3} \\ q_{1} = c_{1} p_{1} + d_{1} p_{2} + e_{1} t_{3} \\ q_{2} = a_{2} p_{1} + b_{2} p_{2} + c_{2} t_{3} \end{array} \\ \hline \end{array} \\[1ex] &\renewcommand{\arraystretch}{1.4} \begin{array}{|c|} \hline 3 \\ \hline \begin{array}{c} t_{1} = a_{1} p_{3} + b_{1} t_{2} \\ r = a_{0} p_{3} + b_{0} t_{2} \end{array} \\ \hline \end{array} \end{align*} } 

\noindent \noindent\textbf{$\bullet$ 19 parameters}:\textbf{ $a_{11},\dots,b_{0}\in\mC$} 

\noindent \noindent\textbf{$\bullet$ Embedding of} $C$: 
\[
C\subset\mP(p_{1},p_{2},t_{3},p_{3},t_{2},u)=\mP(2^{3},3^{2},5).
\]

\noindent \noindent$\bullet$\textbf{ Sing~$T$}: $\mP(p_{1},p_{2},t_{3})\cap T$
consists of four points with $p_{1}\not=0$, say, $\mathsf{p}_{2},\mathsf{q}_{2},\mathsf{r}_{2},\mathsf{s}_{2}$.
$T$ has $1/2(1,1)$-singularities at $\mathsf{p}_{2},\mathsf{q}_{2}$,
$\mathsf{r}_{2},\mathsf{s}_{2}$ and a $1/5(2,3)$-singularity at
the $u$-point.

\noindent \vspace{5pt}

\noindent %
\fbox{\textbf{No.$\,$6865}}

\noindent \vspace{3pt}

\noindent \noindent\textbf{$\bullet$ Sections for $C$}: { \begin{align*} &\renewcommand{\arraystretch}{1.4} \begin{array}{|c|c|c|} \hline \text{weight} & 1 & 2 \\ \hline \text{section} & q_{3} = t_{3} = s_{13} = s_{23} = 0 & \begin{array}{c} \ \ \, s_{11} = a_{11} p_{1} + b_{11} p_{2} + c_{11} p_{4} \\ \ \ \, s_{12} = a_{12} p_{1} + b_{12} p_{2} + c_{12} p_{4} \\ \ \ \, s_{22} = a_{22} p_{1} + b_{22} p_{2} + c_{22} p_{4} \\ t_{1} = a_{1} p_{1} + b_{1} p_{2} + c_{1} p_{4} \\ t_{2} = a_{2} p_{1} + b_{2} p_{2} + c_{2} p_{4} \\ q_{1} = d_{1} p_{1} + e_{1} p_{2} + f_{1} p_{4} \\ q_{2} = d_{2} p_{1} + e_{2} p_{2} + f_{2} p_{4} \end{array} \\ \hline \end{array} \end{align*} } 

\noindent \noindent\textbf{$\bullet$ 21 parameters}:\textbf{ $a_{11},\dots,f_{2}\in\mC$} 

\noindent \noindent\textbf{$\bullet$ Embedding of $C$}: 
\[
C\subset\mP(p_{1},p_{2},p_{4},p_{3},r,u)=\mP(2^{3},3^{2},4).
\]

\noindent \noindent$\bullet$\textbf{ Sing~$T$}: $\mP(p_{1},p_{2},p_{4},u)\cap T$
consists of the $u$-point, and other five points with $p_{1}\not=0$,
say, $\mathsf{p}_{2},\mathsf{q}_{2},\mathsf{r}_{2},\mathsf{s}_{2},\mathsf{t}_{2}$.
$T$ has $1/2(1,1)$-singularities at $\mathsf{p}_{2},\mathsf{q}_{2},\mathsf{r}_{2},\mathsf{s}_{2},\mathsf{t}_{2}$,
and a $1/4(1,3)$-singularity at the $u$-point.

\subsection{Case $h^{0}(\sO_{\mP_{X}}(1))\protect\geq3$\label{subsec:Caseh0ge3}}

In this case, we can easily check that $\text{{\rm Bs}}|\sO_{\mP_{\Sigma}}(b_{i})|_{|T}$
consists of only the $u$-point for any $i$. By this fact and \eqref{eq:SingXT},
$T$ is nonsingular outside of the $u$-point. By LPC, we can determine
the singularity of $T$ at the $u$-point as desired (cf.$\,$Example
\ref{exa:LPC at u-pt}). Thus Claims (A) and (B) follows in this case. 

We can also check that $\text{{\rm Bs}}|\sO_{\mP_{\Sigma}}(a_{i})|_{|X}$
consists of the $u$-point only for any $i$. By the argument in Subsection
\ref{subsec:(A) and (B) =00003D>(C)}, this implies \eqref{eq:0dimT1}.
Thus Claim (C) follows.

\section{\textbf{Proof of Theorem \ref{thm:main1} for Table}~\textbf{\ref{Table3}\label{sec:Proof-of-theorem2}}}

\noindent \indent We note that, in the case where the key variety
of $T$ is the weighted cone $\Pi_{\mathbb{P}}^{14}$ over $\Pi_{\mathbb{P}}^{13}$
(No.$\,$308, 501, 512, 550), the weight one coordinate is unique.
Therefore $T$ is actually a weighted complete intersection of $\Pi_{\mathbb{P}}^{13}$.
Therefore we may assume that $\mathfrak{K}=$$\Pi_{\mP}^{13}$ for
$T$ in this section. 

\noindent \indent For the classes except No.$\,$577, we can argue
similarly to the case that $\mathfrak{K}=\Sigma_{\mP}^{12}$ and $h^{0}(\sO_{X}(1))=1$;
we can show Claims (A) and (B) separately on two charts and the complement
of their union. Claim (C) can be also shown similarly. One difference
is that the choice of the two charts depends on each case. We record
below the choice of the charts in each case. 

\noindent \indent As for No.$\,$577, although the method itself
is straightforward, the computations are much more involved. For details,
we refer to the discussion of No.$\,$577 below, along with the associated
\textit{Mathematica} code \cite{Tak8}.

\noindent \vspace{5pt}

\noindent \textbf{}%
\fbox{\textbf{No.$\,$308}} 

\noindent \vspace{3pt}

\noindent \noindent\textbf{$\bullet$ Sections for $T$} { \begin{align*} &\renewcommand{\arraystretch}{1.4} \begin{array}{|c|c|c|c|c|} \hline \text{weight} & 2 & 3 & 4 & 5 \\ \hline \text{section} & t_{1} = t_{245} = 0 & t_{2} = t_{126} = 0 & p_{2} = t_{124} = t_{136} = 0 & t_{125} = a_{125} p_{1} \\ \hline \end{array} \\[1ex] &\renewcommand{\arraystretch}{1.4} \begin{array}{|c|c|} \hline 6 & 8  \\ \hline \begin{array}{c} t_{123} = a_{123} p_{4} + b_{123} u_{1} \\ t_{135} = a_{135} p_{4} + b_{135} u_{1} \end{array} & s_{1} = a_{1} u_{2}  \\ \hline \end{array} \end{align*} }

\noindent \noindent\textbf{$\bullet$ 6 parameters }$a_{125},\dots,a_{1}\in\mathbb{C}$

\noindent \noindent\textbf{$\bullet$ Embedding of $T$}:
\[
T\subset\mP(p_{1},p_{4},u_{1},p_{3},u_{2},s_{2},s_{3})=\mP(5,6^{2},7,8,9,10).
\]

\noindent \noindent\textbf{$\bullet$ Sing $T$}:\textbf{ }$T\cap\mP(p_{4},u_{1})$
consists of two points, say, $\mathsf{p}_{6}$ and $\mathsf{q_{6}}.$
$T$ has a $1/2(1,1)$-singularity at the point $\mathsf{p}_{2}:=T\cap\mP(p_{4},u_{1},u_{2},s_{3})\cap\{u_{2}\not=0\}$,
a $1/3(2,1)$-singularity at the point $\mathsf{p_{3}:=}T\cap\mP(p_{4},s_{2})$,
a $1/5(2,3)$-singularity at the $p_{1}$-point, and $1/6(1,5)$-singularities
at $\mathsf{p}_{6}$ and $\mathsf{q_{6}}.$ 

\noindent \noindent\textbf{$\bullet$ Two charts: }the $u_{2}$-
and $p_{1}$-charts.

\noindent \vspace{5pt}

\noindent \textbf{}%
\fbox{\textbf{No.$\,$501}} 

\noindent \vspace{3pt}

\noindent \noindent\textbf{$\bullet$ Sections for $T$} { \begin{align*} &\renewcommand{\arraystretch}{1.4} \begin{array}{|c|c|c|c|c|} \hline \text{weight} & 2 & 3 & 4 & 5 \\ \hline \text{section} & t_{1} = t_{245} = 0 & t_{2} = a_2 t_{126} & p_{2} = t_{124} = t_{136} = 0 & p_{1} = t_{125} = 0 \\ \hline \end{array} \\[1ex] &\renewcommand{\arraystretch}{1.4} \begin{array}{|c|} \hline 6 \\ \hline \begin{array}{c} p_{4} = b_{4} u_{1} + c_{4} t_{126}^{2} \\ \quad t_{123} = b_{123} u_{1} + c_{123} t_{126}^{2} \\ \quad t_{135} = b_{135} u_{1} + c_{135} t_{126}^{2} \end{array} \\ \hline \end{array} \end{align*} } 

\noindent \noindent\textbf{$\bullet$ 7 parameters }$a_{2},\dots,c_{135}\in\mathbb{C}$

\noindent \noindent\textbf{$\bullet$ Embedding of $T$}:
\[
T\subset\mP(t_{126},u_{1},p_{3},u_{2},s_{1},s_{2},s_{3})=\mP(3,6,7,8^{2},9,10).
\]

\noindent \noindent\textbf{$\bullet$ Sing $T$}: $T\cap\mP(t_{126},u_{1},s_{2})$
consists of four points, say, $\mathsf{p}_{3},\mathsf{q}_{3},\mathsf{r}_{3},\mathsf{s}_{3}$.
$T$ has a $1/2(1,1)$-singularity at the point $\mathsf{p}_{2}:=T\cap\mP(u_{1},u_{2},s_{1},s_{3})\cap\{u_{1}\not=0\}$,
$1/3(1,2)$-singularities at $\mathsf{p}_{3},\mathsf{q}_{3},\mathsf{r}_{3},\mathsf{s}_{3}$,
and a $1/8(7,1)$-singularity at the $s_{1}$-point.

\noindent \noindent\textbf{$\bullet$ Two charts: }the $u_{1}$-
and $t_{126}$-charts. 

\noindent \vspace{5pt}

\noindent \textbf{}%
\fbox{\textbf{No.$\,$512}} 

\noindent \vspace{3pt}

\noindent \noindent\textbf{$\bullet$ Sections for $T$} { {\small{}\begin{align*} &\renewcommand{\arraystretch}{1.4} \begin{array}{|c|c|c|c|c|} \hline \text{weight} & 2 & 3 & 4 & 5 \\ \hline \text{section} & t_{1} = t_{245} = 0 & \begin{array}{c} p_{2} = a_{2} t_{126} \\ t_{2} = b_2 t_{126} \end{array} & p_{1} = t_{124} = t_{136} = 0 & \begin{array}{c} p_{4} = a_{4} u_{1} \\ t_{125} = a_{125} u_{1} \end{array} \\ \hline \end{array} \\[1ex] &\renewcommand{\arraystretch}{1.4} \begin{array}{|c|} \hline 6 \\ \hline \begin{array}{c} t_{123} = a_{123} p_{3} + b_{123} t_{126}^{2} \\ t_{135} = a_{135} p_{3} + b_{135} t_{126}^{2} \end{array} \\ \hline \end{array} \end{align*} }}

\noindent \noindent\textbf{$\bullet$ 8 parameters }$a_{2},\dots,b_{135}\in\mathbb{C}$

\noindent \noindent\textbf{$\bullet$ Embedding of $T$}:
\[
T\subset\mP(t_{126},u_{1},p_{3},u_{2},s_{1},s_{2},s_{3})=\mP(3,5,6,7^{2},8,9).
\]

\noindent \noindent\textbf{$\bullet$ Sing $T$}: $T\cap\mP(t_{126},p_{3},s_{3})$
consists of three points, say, $\mathsf{p}_{3},\mathsf{q}_{3},\mathsf{r}_{3}$.
$T$ has three $1/3(1,2)$-singularities at $\mathsf{p}_{3},\mathsf{q}_{3},\mathsf{r}_{3}$,
and a $1/5(2,3)$-singularity at the $u_{1}$-point and a $1/7(1,6)$-singularity
at the $s_{1}$-point. 

\noindent \noindent\textbf{$\bullet$ Two charts: }the $t_{126}$-
and $u_{1}$-charts.

\noindent \vspace{5pt}

\noindent \textbf{}%
\fbox{\textbf{No.$\,$550}} 

\noindent \vspace{3pt}

\noindent \noindent\textbf{$\bullet$ Sections for $T$} { {\small{}\begin{align*} &\renewcommand{\arraystretch}{1.4} \begin{array}{|c|c|c|c|c|} \hline \text{weight} & 2 & 3 & 4 & 5 \\ \hline \text{section} & p_{2} = t_{1} = t_{245} = 0 & \begin{array}{c} t_{126} = a_{126} p_1 \\ t_{2} = a_2 p_1\end{array} & \begin{array}{c} p_{4} = a_{4} u_{1} \\ t_{124} = a_{124} u_{1} \\ t_{136} = a_{136} u_{1} \end{array} & t_{125} = a_{125} p_{3} \\ \hline \end{array} \\[1ex] &\renewcommand{\arraystretch}{1.4} \begin{array}{|c|} \hline 6 \\ \hline \begin{array}{c} t_{123} = a_{123} u_{2} + b_{123} s_{1} + c_{123} p_1^{2} \\ t_{135} = a_{135} u_{2} + b_{135} s_{1} + c_{135} p_1^{2} \end{array} \\ \hline \end{array} \end{align*} }} 

\noindent \noindent\textbf{$\bullet$ 12 parameters }$a_{126},\dots,c_{135}\in\mathbb{C}$

\noindent \noindent\textbf{$\bullet$ Embedding of $T$}:
\[
T\subset\mP(p_{1},u_{1},p_{3},u_{2},s_{1},s_{2},s_{3})=\mP(3,4,5,6^{2},7,8).
\]

\noindent \noindent\textbf{$\bullet$ Sing $T$}: $T\cap\mP(p_{1},u_{2},s_{1})\cap\{p_{1}\not=0\}$
consists of three points, say, $\mathsf{p}_{3},\mathsf{q}_{3},\mathsf{r}_{3}$.
$T$ has a $1/2(1,1)$-singularity at the point $\mathsf{p}_{2}:=T\cap\mP(u_{1},u_{2},s_{3})\cap\{u_{1}\not=0\}$,
$1/3(1,2)$-singularities at $\mathsf{p}_{3},\mathsf{q}_{3},\mathsf{r}_{3}$,
a $1/4(1,3)$-singularity at the point $\mathsf{p}_{4}:=T\cap\mP(u_{1},s_{3})$,
and $1/6(1,5)$-singularity at the $s_{1}$-point. 

\noindent \noindent\textbf{$\bullet$ Two charts: }the $p_{1}$-
and $u_{1}$-charts.

\noindent \vspace{5pt}

\noindent \textbf{}%
\fbox{\textbf{No.$\,$872}}\textbf{ }

\noindent \vspace{3pt}

\noindent \noindent\textbf{$\bullet$ Sections for $T$} { 

\noindent {\small{}\begin{align*} &\renewcommand{\arraystretch}{1.4} \begin{array}{|c|c|c|c|c|} \hline \text{weight} & 1 & 2 & 3 & 4 \\ \hline \text{section} & p_{2} = 0 & p_{1} = t_{1} = t_{245} = 0 & \begin{array}{c} t_2=a_2 p_4+b_2 u_1\\ t_{126} =a_{126} p_4+b_{126} u_1 \end{array} & \begin{array}{c} t_{124} = a_{124} p_{3} \\ t_{136} = a_{136} p_{3} \end{array} \\ \hline \end{array} \\[1ex] &\renewcommand{\arraystretch}{1.4} \begin{array}{|c|c|} \hline 5 & 6 \\ \hline t_{125} = a_{125} u_{2} + b_{125} s_{1} & \begin{array}{c} t_{123} = a_{123} s_{2} +b_{123}p_4^2+c_{123}p_4 u_1+d_{123}u_1^2 \\ t_{135} = a_{135} s_{2} +b_{135}p_4^2+c_{135}p_4 u_1+d_{135}u_1^2  \end{array} \\ \hline \end{array} \end{align*} 
}} 

\noindent \noindent\textbf{$\bullet$ 16 parameters }$a_{2},\dots,d_{135}\in\mathbb{C}$

\noindent \noindent\textbf{$\bullet$ Embedding of $T$}:
\[
T\subset\mP(p_{4},u_{1},p_{3},u_{2},s_{1},s_{2},s_{3})=\mP(3^{2},4,5^{2},6,7).
\]

\noindent \noindent\textbf{$\bullet$ Sing $T$}: $T\cap\mP(p_{4},u_{1},s_{2})$
consists of five points $\mathsf{p}_{3},\mathsf{q}_{3},\mathsf{r}_{3},\mathsf{s}_{3},\mathsf{t}_{3}$.
$T$ has $1/3(1,2)$-singularities at $\mathsf{p}_{3},\mathsf{q}_{3},\mathsf{r}_{3},\mathsf{s}_{3},\mathsf{t}_{3}$,
and a $1/5(1,4)$-singularity at the $s_{1}$-point.

\noindent \noindent\textbf{$\bullet$ Two charts:} the $p_{4}$-
and the $u_{1}$-charts.

\noindent \vspace{5pt}

\noindent \textbf{}%
\fbox{\textbf{No.$\,$577}} 

\noindent \vspace{3pt}

\noindent \noindent\textbf{$\bullet$ Sections for $T$} {\begin{align*} &\renewcommand{\arraystretch}{1.4} \begin{array}{|c|c|c|c|c|} \hline \text{weight} & 1 & 2 & 3 & 4 \\ \hline \text{section} & t_{245} = 0 & t_{1} = t_{126} = t_{136} = 0 & \begin{array}{c}t_{2} = b_2 p_2 \\  \,  t_{124} = a_{124}p_2  \\ \, t_{125} = a_{125} p_2\\ \, t_{135}=a_{135} p_2 \end{array} & \begin{array}{c}\, t_{123} = a_{123}p_1  \\ p_{4} = a_{4} p_1 \end{array} \\ \hline \end{array} \\[1ex] &\renewcommand{\arraystretch}{1.4} \begin{array}{|c|} \hline 5 \\ \hline p_{3} = a_{3} u_{1} + b_{3} s_{1} \\ \hline \end{array} \end{align*} 
} 

\noindent \noindent$\bullet$ \textbf{8 parameters $b_{2},\dots,b_{3}\in\mC$}

\noindent \noindent$\bullet$\textbf{ Embedding of $T$}:
\[
T\subset\mP(p_{2},p_{1},u_{1},s_{1},u_{2},s_{2},s_{3})=\mP(3,4,5^{2},6^{2},7).
\]

\noindent \noindent\textbf{$\bullet$ Sing $T$}: $T\cap\mP(p_{2},u_{2},s_{2})$
consists of three points, say, $\mathsf{p}_{3},\mathsf{q}_{3},\mathsf{r}_{3}$,
and $T\cap\mP(u_{1},s_{1})$ consists of two points, say, $\mathsf{p}_{5},\mathsf{q}_{5}$.
$T$ has a $1/2(1,1)$-singularity at the point $\mathsf{p}_{2}:=T\cap\mP(p_{1},u_{2},s_{2})$,
$1/3(1,2)$-singularities at $\mathsf{p}_{3},\mathsf{q}_{3},\mathsf{r}_{3}$,
and $1/5(1,4)$-singularities at $\mathsf{p}_{5},\mathsf{q}_{5}$. 

\noindent \vspace{5pt}

\noindent The full proofs of Claims (A), (B) and (C) for No.$\,$577
is given in the \textit{Mathematica} code \cite{Tak8}. Here we outline
some parts of the proofs. 

\noindent \vspace{3pt}

\noindent \textit{\noindent The smoothness of $T_{\mA}^{o}$}: Let
$W_{\mA}$ be the variety obtained from $T_{\mA}$ by removing the
sections of weights $3$, $4$, $5$. We have the following claim:
\begin{claim}
\label{claim:charts of W}The following open subsets of the $p_{1}$-,
$u_{1}$-, and $s_{1}$-charts are isomorphic to open subsets of the
affine 9-space: 
\begin{equation}
\begin{cases}
\text{the \ensuremath{u_{1}}-chart \ensuremath{\cap}\ensuremath{\{1-p\ensuremath{{}_{1}^{3}t_{2}}-3\ensuremath{p_{1}p_{2}t_{2}}-p\ensuremath{{}_{2}^{3}}t\ensuremath{{}_{2}^{2}\not}=0\},} }\\
\text{the \ensuremath{s_{1}}-chart \ensuremath{\cap\{1+3p_{2}p_{4}t_{135}+p{}_{4}^{3}t_{135}-p{}_{2}^{3}t_{135}^{2}\not=0\},}}\\
\text{the \ensuremath{p_{1}}-chart \ensuremath{\cap}\{\ensuremath{1}+\ensuremath{p_{2}u_{2}\not}=0\}. }
\end{cases}\label{eq:577Wchart}
\end{equation}
\end{claim}
Let $U_{\mA}$ (resp.$\,V_{\mA}$) be the variety obtained from $T_{\mA}$
by removing the section of weight $5$ (resp.~the two sections of
weight $4$). We observe that the restrictions to $U_{\mA}$ (resp.~$V_{\mA}$)
of the open subsets of the $u_{1}$- and the $s_{1}$-charts of $W_{\mA}$
(resp.~the $p_{1}$-chart) as in \eqref{eq:577Wchart} are hypersurfaces.
By the Jacobian criterion, we can check that the hypersurfaces for
the $u_{1}$- and $s_{1}$-charts has at most isolated singularities
(note that, by Claim \ref{claim:charts of W}, we have only to check
the singularities of the hypersurfaces along ${\rm Bs}|\sO(3)|\cup{\rm Bs}|\sO(4)|$).
Moreover, since $w(u_{1})=w(s_{1})=5$, ${\rm Bs}|O(5)|$ is disjoint
from the $u_{1}$- and $s_{1}$-charts. Therefore, since $T_{\mA}$
is obtained from $U_{\mA}$ by cutting a general section of weight
$5$, the $u_{1}$- and $s_{1}$-charts of $T_{\mA}$ are smooth on
the restrictions of the open subsets as in \eqref{eq:577Wchart}.
Similarly, studying the singularities of the hypersurface for the
$p_{1}$-chart, we can check that the $p_{1}$-chart of $T_{\mA}$
is smooth on the restriction of the open subset as in \eqref{eq:577Wchart}.

Now it remains to check that $T_{\mA}^{o}$ is smooth along the locus
$\Delta$ defined by 

\[
\begin{cases}
u_{1}(b_{2}p{}_{1}^{3}p_{2}+b{}_{2}^{2}p{}_{2}^{5}+3b_{2}p_{1}p{}_{2}^{2}u_{1}-u{}_{1}^{3})=0,\\
s_{1}(a_{135}a_{4}^{3}p_{1}^{3}p_{2}-a_{135}^{2}p_{2}^{5}+3a_{135}a_{4}p_{1}p_{2}^{2}s_{1}+s_{1}^{3})=0,\\
p_{1}(p{}_{1}^{2}+p_{2}u_{1})=0.
\end{cases}
\]
By somewhat lengthy calculations, we see that this locus gives only
the locus $T\cap\{u_{1}=s_{1}=p_{1}=0\}$ on $T$, and consists of
three points. Performing LPC there, we see that $T$ has $1/3(1,2)$-singularities
at these three points. Therefore we conclude that $T_{\mA}^{o}$ is
smooth along $\Delta$.

\vspace{3pt}

\noindent \textit{Singularities of $T\cap\{\mathsf{b}=0\}$}: We
have only to show that its intersection with the locus $\Delta'$
defined by 

\[
\begin{cases}
u_{1}(b_{2}p{}_{1}^{3}p_{2}+b{}_{2}^{2}p{}_{2}^{5}+3b_{2}p_{1}p{}_{2}^{2}u_{1}-u{}_{1}^{3})=0, & \text{}\\
p_{1}(p{}_{1}^{2}+p_{2}u_{1})=0
\end{cases}
\]
is a finite set, and $T\cap\{\mathsf{b}=0\}$ has only isolated singularities
on the $p_{1}$- and $u_{1}$-charts outside $\Delta'$. This can
be done in a standard way.

\noindent \vspace{5pt}

\noindent \textbf{}%
\fbox{\textbf{No.$\,$878}} 

\noindent \vspace{3pt}

\noindent \noindent\textbf{$\bullet$ Sections for $T$} { \begin{align*} &\renewcommand{\arraystretch}{1.4} \begin{array}{|c|c|c|} \hline \text{weight} & 1 & 2 \\ \hline \text{section} & t_{245} = 0 & p_{2} = t_{1} = t_{126} = t_{136} = 0 \\ \hline \end{array} \\[1ex] &\renewcommand{\arraystretch}{1.4} \begin{array}{|c|c|} \hline 3 & 4 \\ \hline \begin{array}{c} t_{125} = a_{125}p_1+ b_{125} p_4 \\ t_{135}= a_{135}p_1 + b_{135} p_4 \\ t_{2} = a_2 p_1 + b_2 p_4 \ \ \ \ \\  t_{124} = a_{124} p_1 + b_{124} p_4 \end{array} & \begin{array}{c} p_{3} = a_{3} s_{1} + b_{3} t_{123} \\ u_{1} = c_1 s_{1} + d_1 t_{123} \end{array} \\ \hline \end{array} \end{align*} }

\noindent \noindent\textbf{$\bullet$ 12 parameters }$a_{1},\dots,d_{1}\in\mathbb{C}$

\noindent \noindent\textbf{$\bullet$ Embedding of $T$}:
\[
T\subset\mP(p_{1},p_{4},s_{1},t_{123},u_{2},s_{2},s_{3})=\mP(3^{2},4^{2},5^{2},6).
\]

\noindent \noindent\textbf{$\bullet$ Sing $T$}: $T\cap\mP(p_{1},p_{4},s_{3})$
consists of four points, say, $\mathsf{p}_{3},\mathsf{q}_{3},\mathsf{r}_{3},\mathsf{s}_{3}$,
and $T\cap\mP(s_{1},t_{123})$ consists of two points, say, $\mathsf{p}_{4},\mathsf{q}_{4}$.
$T$ has $1/3(1,2)$-singularities at $\mathsf{p}_{3},\mathsf{q}_{3},\mathsf{r}_{3},\mathsf{s}_{3}$,
and $1/4(1,3)$-singularities at $\mathsf{p}_{4},\mathsf{q}_{4}$.

\noindent \noindent\textbf{$\bullet$ Two charts: }the $p_{1}$-
and the $p_{4}$-charts.

\noindent \vspace{5pt}

\noindent \textbf{}%
\fbox{\textbf{No.$\,$1766}} 

\noindent \vspace{3pt}

\noindent \noindent\textbf{$\bullet$ Sections for $T$} { {\small{}\begin{align*} &\renewcommand{\arraystretch}{1.4} \begin{array}{|c|c|c|} \hline \text{weight} & 1 & 2 \\ \hline \text{section} & p_{2} = t_{245} = 0 & \begin{array}{c} p_{1} = a_{1} t_{136} \\ p_{4} = a_{4} t_{136} \\ t_{1} = d_1 t_{136} \\ t_{126} = a_{126} t_{136} \end{array} \\ \hline \end{array} \\[1ex] &\renewcommand{\arraystretch}{1.4} \begin{array}{|c|c|} \hline 3 & 4 \\ \hline \begin{array}{c} p_{3} = a_{3} u_{1} + b_{3} s_{1} + c_{3} t_{135} \\ t_{2} = d_2 u_{1} +e_2 s_{1} + f_2 t_{135} \\ t_{124} = a_{124} u_{1} + b_{124} s_{1} + c_{124} t_{135} \\ t_{125} = a_{125} u_{1} + b_{125} s_{1} + c_{125} t_{135} \end{array} & s_{2} = a_{2} u_{2} + b_{2} t_{123} + c_{2} t_{136}^{2} \\ \hline \end{array} \end{align*} }} 

\noindent \noindent\textbf{$\bullet$ 19 parameters }$a_{1},\dots,c_{2}\in\mathbb{C}$

\noindent \noindent\textbf{$\bullet$ Embedding of $T$}:
\[
T\subset\mP(t_{136},u_{1},s_{1},t_{135},u_{2},t_{123},s_{3})=\mP(2,3^{3},4^{2},5).
\]

\noindent \noindent\textbf{$\bullet$ Sing $T$}: $T\cap\mP(t_{136},u_{2},t_{123})$
consists of two points, say, $\mathsf{p}_{2},\mathsf{q}_{2}$ with
$t_{136}\not=0$, and $T\cap\mP(u_{1},s_{1},t_{135})$ consists of
five points, say, $\mathsf{p}_{3},\mathsf{q}_{3},\mathsf{r}_{3},\mathsf{s}_{3},\mathsf{t}_{3}$.
$T$ has $1/2(1,1)$-singularities at $\mathsf{p}_{2},\mathsf{q}_{2}$,
and $1/3(1,2)$-singularities at $\mathsf{p}_{3},\mathsf{q}_{3},\mathsf{r}_{3},\mathsf{s}_{3},\mathsf{t}_{3}$.

\noindent \noindent\textbf{$\bullet$ Two charts: }the $u_{1}$-
and the $t_{136}$-charts.

\newpage
\section{\textbf{Tables\label{sec:Tables}}}

\renewcommand{\arraystretch}{2}
{\small
{{ \begin{table}[H]  \rowcolors{2}{}{yellow!20} \caption{Weights of coordinates of $\Sigma^{12}_{\mP}$} \label{Table 1} \begin{tabular}{|c|c|c|c|c|c|c|c|} \hline No.&  $w(\bm{p})$ &  $w(p_4)$ &  $w(\bm{q})$ & $w(r)$ & $w(u)$ & $w(S)$ & $w(\bm{t})$ \\ \hline \hline 360 & $\begin{spmatrix} 7 \\ 8 \\ 9  \end{spmatrix}$ &  $6$ &  $\begin{spmatrix} 6 \\5  \\ 4  \end{spmatrix}$ & $7$ & $8$ & $\begin{spmatrix} 4 & 3 & 2 \\          &2  &1 \\          &           & 0  \end{spmatrix}$ & $\begin{spmatrix}  5 \\4  \\3  \end{spmatrix}$\\\hline 393 & $\begin{spmatrix} 6 \\ 7 \\ 8  \end{spmatrix}$ &  $5$ &  $\begin{spmatrix} 6 \\5  \\ 4  \end{spmatrix}$ & $7$ & $9$ & $\begin{spmatrix} 4 & 3 & 2 \\          &2  &1 \\          &           & 0  \end{spmatrix}$ & $\begin{spmatrix}  5 \\4  \\3  \end{spmatrix}$\\\hline 569 & $\begin{spmatrix} 5 \\ 6 \\ 7  \end{spmatrix}$ &  $3$ &  $\begin{spmatrix} 5 \\4  \\ 3  \end{spmatrix}$ & $6$ & $9$ & $\begin{spmatrix} 4 & 3 & 2 \\          &2  &1 \\          &           & 0  \end{spmatrix}$ & $\begin{spmatrix}  6 \\5  \\4  \end{spmatrix}$\\\hline 574& $\begin{spmatrix} 5 \\6  \\ 7  \end{spmatrix}$ &  $5$ &  $\begin{spmatrix}  5 \\ 4 \\3  \end{spmatrix}$ & $6$ & $7$ & $\begin{spmatrix} 4 & 3 & 2\\         & 2 &1 \\          &           &0  \end{spmatrix}$ & $\begin{spmatrix}  4\\ 3 \\ 2  \end{spmatrix}$\\\hline 642 & $\begin{spmatrix} 4 \\5  \\6  \end{spmatrix}$ &  $3$ &  $\begin{spmatrix} 6 \\ 5 \\ 4 \end{spmatrix}$ & $7$ & $11$ & $\begin{spmatrix} 4 & 3 & 2 \\          & 2 &1 \\          &   & 0 \end{spmatrix}$ & $\begin{spmatrix} 5 \\4  \\3  \end{spmatrix}$\\\hline 644 & $\begin{spmatrix} 4 \\5  \\6  \end{spmatrix}$ &  $4$ &  $\begin{spmatrix} 6 \\5  \\4  \end{spmatrix}$ & $7$ & $10$ & $\begin{spmatrix} 4 & 3 & 2\\          & 2 &1 \\          &  & 0 \end{spmatrix}$ & $\begin{spmatrix} 4 \\3  \\2  \end{spmatrix}$\\\hline 1091 & $\begin{spmatrix} 6 \\ 7 \\ 8  \end{spmatrix}$ &  $5$ &  $\begin{spmatrix} 6  \\5  \\ 4  \end{spmatrix}$ & $7$ & $9$ & $\begin{spmatrix} 4 & 3 & 2\\          & 2 & 1\\          &  &0  \end{spmatrix}$ & $\begin{spmatrix} 5 \\4  \\3  \end{spmatrix}$\\\hline 1181 & $\begin{spmatrix} 3 \\4  \\5  \end{spmatrix}$ &  $2$ &  $\begin{spmatrix} 6 \\5  \\4  \end{spmatrix}$ & $7$ & $12$ & $\begin{spmatrix} 4& 3 & 2\\          & 2 & 1\\          &  & 0 \end{spmatrix}$ & $\begin{spmatrix} 5 \\4  \\3  \end{spmatrix}$\\\hline 1185 & $\begin{spmatrix} 4 \\5  \\6  \end{spmatrix}$ &  $2$ &  $\begin{spmatrix} 4 \\ 3 \\2  \end{spmatrix}$ & $5$ & $8$ & $\begin{spmatrix} 4 & 3 & 2 \\          & 2 & 1 \\          &  &0  \end{spmatrix}$ & $\begin{spmatrix} 6 \\5  \\4  \end{spmatrix}$\\\hline
1186 & $\begin{spmatrix} 4 \\5  \\6  \end{spmatrix}$ &  $3$ &  $\begin{spmatrix} 4 \\3  \\2  \end{spmatrix}$ & $5$ & $7$ & $\begin{spmatrix} 4& 3 & 2\\          & 2 & 1\\          &  & 0 \end{spmatrix}$ & $\begin{spmatrix} 5 \\4  \\3  \end{spmatrix}$\\\hline 1218 & $\begin{spmatrix} 4 \\ 5 \\6  \end{spmatrix}$ &  $5$ &  $\begin{spmatrix} 4 \\3  \\2  \end{spmatrix}$ & $5$ & $5$ & $\begin{spmatrix} 4& 3 & 2\\          & 2 & 1\\          &           &0  \end{spmatrix}$ & $\begin{spmatrix} 3 \\2  \\1  \end{spmatrix}$\\\hline 1253 & $\begin{spmatrix} 4 \\4  \\5  \end{spmatrix}$ &  $3$ &  $\begin{spmatrix}  4\\ 4  \\ 3  \end{spmatrix}$ & $5$ & $7$ & $\begin{spmatrix} 2& 2 & 1\\         & 2 & 1\\          &           &0  \end{spmatrix}$ & $\begin{spmatrix} 3 \\3  \\2  \end{spmatrix}$\\\hline 1413 & $\begin{spmatrix} 4 \\4  \\5  \end{spmatrix}$ &  $3$ &  $\begin{spmatrix} 3 \\3  \\2  \end{spmatrix}$ & $4$& $5$ & $\begin{spmatrix} 2& 2 & 1\\          &2  &1 \\          &           &0  \end{spmatrix}$ & $\begin{spmatrix}  3\\3  \\2  \end{spmatrix}$\\\hline 2422 & $\begin{spmatrix} 2 \\3  \\4  \end{spmatrix}$ &  $3$ &  $\begin{spmatrix} 4 \\3  \\2  \end{spmatrix}$ & $5$ & $7$ & $\begin{spmatrix} 4& 3 & 2\\          &2  & 1\\          &           &0  \end{spmatrix}$ & $\begin{spmatrix} 3 \\2  \\1  \end{spmatrix}$\\\hline 4850& $\begin{spmatrix} 4 \\5  \\6  \end{spmatrix}$ &  $1$ &  $\begin{spmatrix} 6 \\5  \\4  \end{spmatrix}$ & $7$ & $13$ & $\begin{spmatrix} 4& 3 & 2\\          & 2 & 1\\          &           &0  \end{spmatrix}$ & $\begin{spmatrix} 7 \\6  \\5  \end{spmatrix}$\\\hline 4938& $\begin{spmatrix} 3 \\4  \\5  \end{spmatrix}$ &  $1$ &  $\begin{spmatrix} 5 \\4  \\3  \end{spmatrix}$ & $6$ & $11$ & $\begin{spmatrix} 4& 3 & 2\\          & 2 & 1\\          &           &0  \end{spmatrix}$ & $\begin{spmatrix} 6 \\5  \\4  \end{spmatrix}$\\\hline 5202& $\begin{spmatrix} 2 \\3  \\4  \end{spmatrix}$ &  $1$ &  $\begin{spmatrix} 4 \\3  \\2  \end{spmatrix}$ & $5$& $9$& $\begin{spmatrix} 4& 3 & 2\\          & 2 & 1\\          &           &0  \end{spmatrix}$ & $\begin{spmatrix} 5 \\4  \\3  \end{spmatrix}$\\\hline 5859& $\begin{spmatrix} 1 \\2  \\3  \end{spmatrix}$ &  $2$ &  $\begin{spmatrix} 4 \\3  \\2  \end{spmatrix}$ & $5$& $8$ & $\begin{spmatrix} 4& 3 & 2\\          & 2 &1 \\          &           &0  \end{spmatrix}$ & $\begin{spmatrix} 3 \\2  \\1  \end{spmatrix}$\\\hline 5866& $\begin{spmatrix} 2 \\2  \\3  \end{spmatrix}$ &  $1$ &  $\begin{spmatrix}  3\\3  \\2  \end{spmatrix}$ & $4$& $7$& $\begin{spmatrix} 2& 2 &1 \\          &2  &1 \\          &           &0  \end{spmatrix}$ & $\begin{spmatrix} 3 \\3  \\2  \end{spmatrix}$\\\hline 6860& $\begin{spmatrix} 2 \\2  \\3  \end{spmatrix}$ &  $1$ &  $\begin{spmatrix}  2\\2  \\1  \end{spmatrix}$ & $3$& $5$& $\begin{spmatrix} 2& 2 & 1\\          &2  &1 \\          &           &0  \end{spmatrix}$ & $\begin{spmatrix} 3 \\3  \\2  \end{spmatrix}$\\\hline 6865& $\begin{spmatrix} 2 \\ 2 \\3  \end{spmatrix}$ &  $2$ &  $\begin{spmatrix} 2 \\2  \\1  \end{spmatrix}$ & $3$& $4$ & $\begin{spmatrix} 2& 2 &1 \\          &2  &1 \\          &           &0  \end{spmatrix}$ & $\begin{spmatrix} 2 \\2  \\1  \end{spmatrix}$\\\hline 11004& $\begin{spmatrix} 1 \\2  \\3  \end{spmatrix}$ &  $1$ &  $\begin{spmatrix} 3 \\2  \\1  \end{spmatrix}$ & $4$ & $7$ & $\begin{spmatrix} 4& 3 & 2\\          & 2 & 1\\          &           &0  \end{spmatrix}$ & $\begin{spmatrix} 4 \\3  \\2  \end{spmatrix}$\\\hline 16227& $\begin{spmatrix} 1 \\1  \\2  \end{spmatrix}$ &  $1$&  $\begin{spmatrix} 2 \\2  \\1  \end{spmatrix}$ & $3$& $5$& $\begin{spmatrix} 2& 2 & 1\\          & 2 & 1\\          &           &0  \end{spmatrix}$ & $\begin{spmatrix} 2 \\2  \\1  \end{spmatrix}$\\\hline \end{tabular} \label{tab:keyvarwt}  \end{table} }}
}\vspace{3pt}
\noindent Notation for Tables ~\ref{tab:keyvarwt}. 

\vspace{3pt}

\noindent $\bullet$ $w(*):=\text{the weight of the coordinate \ensuremath{*}}$.

\vspace{3pt}

\noindent $\bullet w(\bm{p}), w(\bm{q}), w(S)$, and $w(\bm{t})$ represent the weights of the entries of $\bm{p}$, $\bm{q}, S$ and $\bm{t}$ introduced in Subsection \ref{subsec:The-key-varietySigma}.
\begin{table}[H]
  \caption{Descriptions of $\mQ$-Fano threefolds $X$ in $\Sigma_{\mP}^{12}$}
  \label{Table 2}
  \centering
  \renewcommand{\arraystretch}{1.2}
  \rowcolors{2}{}{yellow!20}
  \resizebox{\textwidth}{!}{
    \begin{tabular}{|c|c|c|c|}
    \hline
    No. & $\mP_X$ & Baskets of singularities & $X\subset\Sigma:=\Sigma_{\mP}^{12}$ \\
    \hline
    \hline
    360 & $\mP(1,4,5,6,7^{2},8,9)$ & $2\times1/4(1,1,3),1/6(1,1,5),1/7(1,2,5)$ & $\Sigma\cap(2)^{2}\cap(3)^{2}\cap(4)^{2}\cap(5)\cap(6)\cap(8)$ \\
    393 & $\mP(1,4,5^{2},6,7,8,9)$ & $1/2(1,1,1),1/5(1,1,4),1/5(1,2,3),1/9(1,4,5)$ & $\Sigma\cap(2)^{2}\cap(3)^{2}\cap(4)^{2}\cap(5)\cap(6)\cap(7)$ \\
    569 & $\mP(1,3,4,5^{2},6,7,9)$ & $2\times1/3(1,1,2),1/5(1,2,3),1/9(1,4,5)$ & $\Sigma\cap(2)^{2}\cap(3)^{2}\cap(4)^{2}\cap(5)\cap(6)^{2}$ \\
    574 & $\mP(1,3,4,5^{2},6,7^{2})$ & $1/3(1,1,2),1/5(1,1,4),1/5(1,2,3),1/7(1,3,4)$ & $\Sigma\cap(2)^{3}\cap(3)^{2}\cap(4)^{2}\cap(5)\cap(6)$ \\
    642 & $\mP(1,3,4^{2},5,6,7,11)$ & $1/2(1,1,1),1/3(1,1,2),1/4(1,1,3),1/11(1,4,7)$ & $\Sigma\cap(2)^{2}\cap(3)^{2}\cap(4)^{2}\cap(5)^{2}\cap(6)$ \\
    644 & $\mP(1,3,4^{2},5,6,7,10)$ & $1/2(1,1,1),2\times1/4(1,1,3),1/10(1,3,7)$ & $\Sigma\cap(2)^{3}\cap(3)\cap(4)^{3}\cap(5)\cap(6)$ \\
    1091 & $\mP(1,2,5,6,7^{2},8,9)$ & $2\times1/2(1,1,1),1/7(1,1,6),1/9(1,2,7)$ & $\Sigma\cap(2)\cap(3)^{2}\cap(4)^{3}\cap(5)^{2}\cap(6)$ \\
    1181 & $\mP(1,2,3,4,5^{2},7,12)$ & $1/2(1,1,1),1/4(1,1,3),1/12(1,5,7)$ & $\Sigma\cap(2)^{2}\cap(3)^{2}\cap(4)^{3}\cap(5)\cap(6)$ \\
    1185 & $\mP(1,2,3,4,5^{2},6,8)$ & $3\times1/2(1,1,1),1/5(1,1,4),1/8(1,3,5)$ & $\Sigma\cap(2)^{3}\cap(3)\cap(4)^{3}\cap(5)\cap(6)$ \\
    1186 & $\mP(1,2,3,4,5^{2},6,7)$ & $2\times1/2(1,1,1),1/3(1,1,2),1/5(1,1,4),1/7(1,2,5)$ & $\Sigma\cap(2)^{2}\cap(3)^{3}\cap(4)^{3}\cap(5)$ \\
    1218 & $\mP(1,2,3,4,5^{3},6)$ & $2\times1/2(1,1,1),1/5(1,1,4),2\times1/5(1,2,3)$ & $\Sigma\cap(1)\cap(2)^{3}\cap(3)^{2}\cap(4)^{2}\cap(5)$ \\
    1253 & $\mP(1,2,3,4^{2},5^{2},7)$ & $2\times1/2(1,1,1),2\times1/4(1,1,3),1/7(1,2,5)$ & $\Sigma\cap(1)\cap(2)^{3}\cap(3)^{3}\cap(4)^{2}$ \\
    1413 & $\mP(1,2,3^{2},4^{2},5^{2})$ & $2\times1/2(1,1,1),2\times1/3(1,1,2),1/4(1,1,3),1/5(1,2,3)$ & $\Sigma\cap(1)\cap(2)^{4}\cap(3)^{3}\cap(4)$ \\
    2422 & $\mP(1,2^{2},3^{2},4,5,7)$ & $4\times1/2(1,1,1),1/3(1,1,2),1/7(1,2,5)$ & $\Sigma\cap(1)\cap(2)^{3}\cap(3)^{3}\cap(4)^{2}$ \\
    4850 & $\mP(1^{2},4,5,6^{2},7,13)$ & $1/2(1,1,1),1/13(1,6,7)$ & $\Sigma\cap(2)^{2}\cap(3)\cap(4)^{2}\cap(5)^{2}\cap(6)\cap(7)$ \\
    4938 & $\mP(1^{2},3,4,5^{2},6,11)$ & $1/3(1,1,2),1/11(1,5,6)$ & $\Sigma\cap(2)^{2}\cap(3)^{2}\cap(4)^{3}\cap(5)\cap(6)$ \\
    5202 & $\mP(1^{2},2,3,4^{2},5,9)$ & $2\times1/2(1,1,1),1/9(1,4,5)$ & $\Sigma\cap(2)^{3}\cap(3)^{3}\cap(4)^{2}\cap(5)$ \\
    5859 & $\mP(1^{2},2^{2},3^{2},5,8)$ & $2\times1/2(1,1,1),1/8(1,3,5)$ & $\Sigma\cap(1)\cap(2)^{4}\cap(3)^{2}\cap(4)^{2}$ \\
    5866 & $\mP(1^{2},2^{2},3^{2},4,7)$ & $3\times1/2(1,1,1),1/7(1,3,4)$ & $\Sigma\cap(1)\cap(2)^{5}\cap(3)^{3}$ \\
    6860 & $\mP(1^{2},2^{3},3^{2},5)$ & $4\times1/2(1,1,1),1/5(1,2,3)$ & $\Sigma\cap(1)^{2}\cap(2)^{5}\cap(3)^{2}$ \\
    6865 & $\mP(1^{2},2^{3},3^{2},4)$ & $5\times1/2(1,1,1),1/4(1,1,3)$ & $\Sigma\cap(1)^{2}\cap(2)^{7}$ \\
    11004 & $\mP(1^{3},2,3^{2},4,7)$ & $1/7(1,3,4)$ & $\Sigma\cap(1)\cap(2)^{4}\cap(3)^{2}\cap(4)^{2}$ \\
    16227 & $\mP(1^{4},2^{2},3,5)$ & $1/5(1,2,3)$ & $\Sigma\cap(1)^{3}\cap(2)^{6}$ \\
    \hline
    \end{tabular}
  }
\end{table}

\begin{table}[H]
  \caption{Weights of coordinates of $\Pi_{\mP}^{13}$}
  \label{tab:Weights-of-coordinatesPi}
  \centering
  \renewcommand{\arraystretch}{1.2}
  \rowcolors{2}{}{yellow!20}
  \resizebox{\textwidth}{!}{
    \begin{tabular}{|c|c|c|}
    \hline
    weights & No.308,~501,~512,~550,~872 & No.577,~ 878,~ 1766 \\
    \hline
    \hline
    $\begin{spmatrix}w(u_1) & w(u_2) \\ w(p_1) & w(p_3) \\ w(p_2) & w(p_4) \end{spmatrix}$ & $\begin{spmatrix}d-1 & d+1 \\ d-2 & d \\ d-3 & d-1 \end{spmatrix}$ & $\begin{spmatrix}d+1 & d+2 \\ d & d+1 \\ d-1 & d \end{spmatrix}$ \\
    \hline
    $\begin{spmatrix}w(s_1) & w(s_2)& w(s_3)\end{spmatrix}$ & \multicolumn{2}{c|}{$\begin{spmatrix}d+1 & d+2 & d+3\end{spmatrix}$} \\
    \hline
    $\begin{spmatrix}w(t_1) & w(t_2)\end{spmatrix}$ & \multicolumn{2}{c|}{$\begin{spmatrix}2 & 3\end{spmatrix}$} \\
    \hline
    $\begin{spmatrix}w(t_{123}) & w(t_{125}) \\ w(t_{124}) & w(t_{126})\end{spmatrix}$ & $\begin{spmatrix}6 & 5 \\ 4 & 3\end{spmatrix}$ & $\begin{spmatrix}4 & 3 \\ 3 & 2\end{spmatrix}$ \\
    \hline
    $\begin{spmatrix}w(t_{135}) & w(t_{136}) & w(t_{245})\end{spmatrix}$ & $\begin{spmatrix}6 & 4 & 2\end{spmatrix}$ & $\begin{spmatrix}3 & 2 & 1\end{spmatrix}$ \\
    \hline
    \end{tabular}
  }
  \vspace{3pt}
  \small
  In Table~\ref{tab:Weights-of-coordinatesPi}, $d=7,7,6,5,4,4,3,2$ for No.~308, 501, 512, 550, 872, 577, 878, 1766, respectively.
  For No.~308, 501, 512, 550, we consider the cone $\Pi_{\mP}^{14}$ over $\Pi_{\mP}^{13}$ with a coordinate of weight 1 being added.
\end{table}

\begin{table}[H]
  \caption{Descriptions of $\mathbb{Q}$-Fano threefolds $X$ in $\Pi_{\mP}^{14}$ or $\Pi_{\mP}^{13}$}
  \label{Table3}
  \centering
  \renewcommand{\arraystretch}{1.2}
  \rowcolors{2}{}{yellow!20}
  \resizebox{\textwidth}{!}{
    \begin{tabular}{|c|c|c|c|}
    \hline
    No. & $\mathbb{P}_{X}$ & Basket of singularities & $X\subset\Pi_{\mP}^{14}$ or $\Pi_{\mP}^{13}$ \\
    \hline
    308 & $\mathbb{P}(1,5,6^{2},7,8,9,10)$ & $1/2(1,1,1),1/3(1,1,2),1/5(1,2,3),2\times1/6(1,1,5)$ & $\Pi_{\mathbb{P}}^{14}\cap(2)^{2}\cap(3)^{2}\cap(4)^{3}\cap(5)\cap(6)^{2}\cap(8)$ \\
    501 & $\mathbb{P}(1,3,6,7,8^{2},9,10)$ & $1/2(1,1,1),4\times1/3(1,1,2),1/8(1,1,7)$ & $\Pi_{\mathbb{P}}^{14}\cap(2)^{2}\cap(3)\cap(4)^{3}\cap(5)^{2}\cap(6)^{3}$ \\
    512 & $\mathbb{P}(1,3,5,6,7^{2},8,9)$ & $3\times1/3(1,1,2),1/5(1,2,3),1/7(1,1,6)$ & $\Pi_{\mathbb{P}}^{14}\cap(2)^{2}\cap(3)^{2}\cap(4)^{3}\cap(5)^{2}\cap(6)^{2}$ \\
    550 & $\mathbb{P}(1,3,4,5,6^{2},7,8)$ & $1/2(1,1,1),3\times1/3(1,1,2),1/4(1,1,3),1/6(1,1,5)$ & $\Pi_{\mathbb{P}}^{14}\cap(2)^{3}\cap(3)^{2}\cap(4)^{3}\cap(5)\cap(6)^{2}$ \\
    872 & $\mathbb{P}(1,3^{2},4,5^{2},6,7)$ & $5\times1/3(1,1,2),1/5(1,1,4)$ & $\Pi_{\mathbb{P}}^{13}\cap(2)^{3}\cap(3)^{2}\cap(4)^{2}\cap(5)\cap(6)^{2}$ \\
    577 & $\mathbb{P}(1,3,4,5^{2},6^{2},7)$ & $1/2(1,1,1),3\times1/3(1,1,2),2\times1/5(1,1,4)$ & $\Pi_{\mathbb{P}}^{13}\cap(2)^{3}\cap(3)^{4}\cap(4)^{2}\cap(5)$ \\
    878 & $\mathbb{P}(1,3^{2},4^{2},5^{2},6)$ & $4\times1/3(1,1,2),2\times1/4(1,1,3)$ & $\Pi_{\mathbb{P}}^{13}\cap(2)^{4}\cap(3)^{4}\cap(4)^{2}$ \\
    1766 & $\mathbb{P}(1,2,3^{3},4^{2},5)$ & $2\times1/2(1,1,1),5\times1/3(1,1,2)$ & $\Pi_{\mathbb{P}}^{13}\cap(1)\cap(2)^{4}\cap(3)^{4}\cap(4)$ \\
    \hline
    \end{tabular}
  }
\end{table}

\appendix

\section{\textbf{Properties of the key varieties\label{sec:Equations-of-the key variety}}}

In this section, we review the definitions and the properties of the
key varieties $\Sigma_{\mA}^{13}$ and $\Pi_{\mA}^{14}$ following
\cite{Tak9} and \cite{Tak6}.

\subsection{Definition of the key variety $\Sigma_{\mA}^{13}$\label{subsec:The-key-varietySigma}}
\begin{defn}
\label{defn:Key} We consider the $18$-dimensional affine space $\mA_{\Sigma}$
with coordinates implemented in the following form with vectors and
matrices: 
\begin{align*}
  & \bm{p} := \begin{pmatrix} p_1 \\ p_2 \\ p_3 \end{pmatrix},\
    p_4,\
    \bm{q} := \begin{pmatrix} q_1 \\ q_2 \\ q_3 \end{pmatrix},\
    r,\, u, \\
  & S := \begin{pmatrix}
        s_{11} & s_{12} & s_{13} \\
              & s_{22} & s_{23} \\
              &        & s_{33}
      \end{pmatrix},\
    \bm{t} := \begin{pmatrix} t_1 \\ t_2 \\ t_3 \end{pmatrix},
\end{align*}
where $S$ is a symmetric matrix. We also introduce the following
anti-symmetric matrix:
\[
A_{\bm{q}}:=\begin{pmatrix}
0 & -q_{3} & q_{2}\\
 & 0 & -q_{1}\\
 &  & 0
\end{pmatrix}.
\]
\begin{enumerate}[(1)]

\item In the affine space $\mA_{\Sigma}$ with the above coordinates,
we define $\Sigma_{\mA}^{14}$ to be the affine scheme defined by
the following nine polynomials $F_{1},\dots,F_{9}$: 

\[
\begin{cases}
F_{1}:={\empty^{t}\bm{p}}{\bm{q}},\\
\empty^{t}\!\begin{pmatrix}
F_{2} & F_{3} & F_{4}\end{pmatrix}:=\left(rI+A_{\bm{q}}S\right)\bm{p}+p_{4}A_{\bm{q}}\bm{t},\\
F_{5}:={\empty^{t}\bm{p}}S\bm{p}+p_{4}{\empty^{t}\bm{p}}\bm{t},\\
\empty^{t}\!\begin{pmatrix}
F_{6} & F_{7} & F_{8}\end{pmatrix}:=u\bm{p}-\left(rI-A_{\bm{q}}S\right)A_{\bm{q}}\bm{t},\\
F_{9}:=up_{4}+\left(r^{2}+{\empty^{t}\bm{q}}S^{\dagger}\bm{q}\right),
\end{cases}
\]
where $I$ is the $3\times3$ identity matrix and $S^{\dagger}$ is
the adjoint matrix of $S$. 

\item We set 
\[
\Sigma_{\mA}^{13}:=\Sigma_{\mA}^{14}\cap\{s_{33}=1\}.
\]
\end{enumerate}
\end{defn}

\subsection{Definition of the key variety $\Pi_{\mA}^{15}${\small{}\label{subsec:Definition-ofPi}}}
\begin{defn}
\label{def:Pi15} Let $\mA_{\Pi}$ be the affine $19$-space with
coordinates 
\begin{equation}
p_{1},\dots,p_{4},u_{1},u_{2},s_{1},s_{2},s_{3},t_{1},t_{2},t_{123},t_{124},t_{125},t_{126},t_{135},t_{136},t_{245},t_{246}.\label{eq:coordpi}
\end{equation}

\begin{enumerate}[(1)]

\item In $\mA_{\Pi}$, we define $\Pi_{\mA}^{15}$ to be the affine
scheme defined by the following nine polynomials $G_{1},\dots,G_{9}$:
\small{\begin{equation*}
\begin{aligned}
G_{1}:={}&\ u_{1}s_{1} - p_{1}s_{2} - p_{2}s_{3} + t_{126}(p_{2}p_{3} - p_{1}p_{4}) \\ &+ t_{136}(-p_{1}^{2} - t_{1}p_{2}^{2} - p_{2}u_{1}) + t_{245}(-2p_{1}p_{3} - 2t_{1}p_{2}p_{4} - p_{4}u_{1} - p_{2}u_{2}) \\ &+ t_{246}(-p_{3}^{2} - t_{1}p_{4}^{2} - p_{4}u_{2}),
\end{aligned}
\end{equation*}
\begin{equation*}
\begin{aligned}
G_{2}:={}&\ u_{2}s_{1} - p_{3}s_{2} - p_{4}s_{3} + t_{125}(-p_{2}p_{3} + p_{1}p_{4}) \\ &+ t_{135}(p_{1}^{2} + t_{1}p_{2}^{2} + p_{2}u_{1}) + t_{136}(2p_{1}p_{3} + 2t_{1}p_{2}p_{4} + p_{4}u_{1} + p_{2}u_{2}) \\ &+ t_{245}(p_{3}^{2} + t_{1}p_{4}^{2} + p_{4}u_{2}), 
\end{aligned}
\end{equation*}
\begin{equation*}
\begin{aligned}
G_{3} := {}& -t_{2}p_{2}s_{1} + (t_{1}p_{2} + u_{1})s_{2} - p_{1}s_{3} + t_{124}(p_{2}p_{3} - p_{1}p_{4}) \\ &+ t_{136}(-p_{2}^{2}t_{2} - p_{1}u_{1}) + t_{245}(-2t_{2}p_{2}p_{4} - p_{3}u_{1} - p_{1}u_{2}) \\ &+ t_{246}(-t_{2}p_{4}^{2} - p_{3}u_{2}),
\end{aligned}
\end{equation*}
\begin{equation*}
\begin{aligned}
G_{4} := {}& -t_{2}p_{4}s_{1} + (t_{1}p_{4} + u_{2})s_{2} - p_{3}s_{3} + t_{123}(-p_{2}p_{3} + p_{1}p_{4}) \\ &+ t_{135}(t_{2}p_{2}^{2} + p_{1}u_{1}) + t_{136}(2t_{2}p_{2}p_{4} + p_{3}u_{1} + p_{1}u_{2}) \\ &+ t_{245}(t_{2}p_{4}^{2} + p_{3}u_{2}), 
\end{aligned}
\end{equation*}
\begin{equation*}
\begin{aligned}
G_{5} :={}& -t_{2}p_{1}s_{1} + (t_{1}p_{1} - t_{2}p_{2})s_{2} + (t_{1}p_{2} + u_{1})s_{3} \\ &+ t_{123}(t_{1}p_{2}^{2} + p_{1}^{2} + p_{2}u_{1}) + t_{124}(t_{1}p_{2}p_{4} + p_{1}p_{3} + p_{2}u_{2}) \\ &- t_{125}(t_{2}p_{2}^{2} + p_{1}u_{1}) - t_{126}(t_{2}p_{2}p_{4} + p_{1}u_{2}) \\ &- t_{136}(t_{1}p_{2}u_{1} - t_{2}p_{1}p_{2} + u_{1}^{2}) \\ &- t_{245}(t_{1}p_{4}u_{1} + t_{1}p_{2}u_{2} - t_{2}p_{2}p_{3} - t_{2}p_{1}p_{4} + 2u_{1}u_{2}) \\ &- t_{246}(t_{1}p_{4}u_{2} - t_{2}p_{3}p_{4} + u_{2}^{2}),
\end{aligned}
\end{equation*}
\begin{equation*}
\begin{aligned}
G_{6} :={}& -t_{2}p_{3}s_{1} + (t_{1}p_{3} - t_{2}p_{4})s_{2} + (t_{1}p_{4} + u_{2})s_{3} \\ &+ t_{123}(p_{1}p_{3} + t_{1}p_{2}p_{4} + p_{4}u_{1}) + t_{124}(p_{3}^{2} + t_{1}p_{4}^{2} + p_{4}u_{2}) \\ &- t_{125}(t_{2}p_{2}p_{4} + p_{3}u_{1}) - t_{126}(t_{2}p_{4}^{2} + p_{3}u_{2}) \\ &+ t_{135}(t_{1}p_{2}u_{1} - t_{2}p_{1}p_{2} + u_{1}^{2}) \\ &+ t_{136}(t_{1}p_{4}u_{1} + t_{1}p_{2}u_{2} - t_{2}p_{2}p_{3} - t_{2}p_{1}p_{4} + 2u_{1}u_{2}) \\ &+ t_{245}(t_{1}p_{4}u_{2} - t_{2}p_{3}p_{4} + u_{2}^{2}),
\end{aligned}
\end{equation*}
\begin{equation*}
\begin{aligned}
G_{7} := {}&-s_{2}^{2} + s_{1}s_{3} + (t_{123}p_{2} + t_{124}p_{4})s_{1} - (t_{125}p_{2} + t_{126}p_{4})s_{2} \\ &+ (t_{136}^{2} - t_{135}t_{245})(p_{1}^{2} - 2p_{2}u_{1} - t_{1}p_{2}^{2}) \\ &+ (t_{245}^{2} - t_{136}t_{246})(p_{3}^{2} - 2p_{4}u_{2} - t_{1}p_{4}^{2}) \\ &+ (t_{135}t_{246} - t_{136}t_{245})(p_{2}u_{2} + p_{4}u_{1} - p_{1}p_{3} + t_{1}p_{2}p_{4}) \\ &+ (t_{123}t_{136} - t_{124}t_{135})p_{2}^{2} + 2(t_{123}t_{245} - t_{124}t_{136})p_{2}p_{4} \\ &+ (t_{123}t_{246} - t_{124}t_{245})p_{4}^{2} + (t_{126}t_{135} - t_{125}t_{136})p_{1}p_{2} \\ &+ (t_{126}t_{136} - t_{125}t_{245})(p_{1}p_{4} + p_{2}p_{3}) + (t_{126}t_{245} - t_{125}t_{246})p_{3}p_{4}, 
\end{aligned}
\end{equation*}
\begin{equation*}
\begin{aligned}
G_{8} :={}&\ s_{2}s_{3} + t_{1}s_{1}s_{2} - t_{2}s_{1}^{2} + (t_{123}p_{1} + t_{124}p_{3})s_{1} - (t_{125}p_{1} + t_{126}p_{3})s_{2} \\ &+ (t_{136}^{2} - t_{135}t_{245})(-2p_{1}u_{1} - 2t_{1}p_{1}p_{2} + t_{2}p_{2}^{2}) \\ &+ (t_{245}^{2} - t_{136}t_{246})(-2p_{3}u_{2} - 2t_{1}p_{3}p_{4} + t_{2}p_{4}^{2}) \\ &+ (t_{135}t_{246} - t_{136}t_{245})(p_{1}u_{2} + p_{3}u_{1} + t_{1}(p_{1}p_{4} + p_{2}p_{3}) - t_{2}p_{2}p_{4}) \\ &+ (t_{126}t_{135} - t_{125}t_{136})p_{1}^{2} + 2(t_{126}t_{136} - t_{125}t_{245})p_{1}p_{3} \\ &+ (t_{126}t_{245} - t_{125}t_{246})p_{3}^{2} + (t_{123}t_{136} - t_{124}t_{135})p_{1}p_{2} \\ &+ (t_{123}t_{245} - t_{124}t_{136})(p_{1}p_{4} + p_{2}p_{3}) + (t_{123}t_{246} - t_{124}t_{245})p_{3}p_{4}, 
\end{aligned}
\end{equation*}
\begin{equation*}
\begin{aligned}
G_{9} :={}&\ s_{3}^{2} + t_{1}s_{2}^{2} - t_{2}s_{1}s_{2} + (t_{123}u_{1} + t_{124}u_{2})s_{1} - (t_{125}u_{1} + t_{126}u_{2})s_{2} \\ &+ (-t_{136}^{2} + t_{135}t_{245})(2t_{2}p_{1}p_{2} + u_{1}^{2}) \\ &+ (-t_{136}t_{245} + t_{135}t_{246})(t_{2}(p_{2}p_{3} + p_{1}p_{4}) + u_{1}u_{2}) \\ &+ (-t_{245}^{2} + t_{136}t_{246})(2t_{2}p_{3}p_{4} + u_{2}^{2}) \\ &- t_{1}\big((t_{123}t_{136} - t_{124}t_{135})p_{2}^{2} + 2(t_{123}t_{245} - t_{124}t_{136})p_{2}p_{4} + (t_{123}t_{246} - t_{124}t_{245})p_{4}^{2}\big) \\ &+ t_{2}\big((t_{125}t_{136} - t_{126}t_{135})p_{2}^{2} + 2(t_{125}t_{245} - t_{126}t_{136})p_{2}p_{4} + (t_{125}t_{246} - t_{126}t_{245})p_{4}^{2}\big) \\ &+ (t_{124}t_{135} - t_{123}t_{136})p_{1}^{2} - 2(-t_{124}t_{136} + t_{123}t_{245})p_{1}p_{3} + (t_{124}t_{245} - t_{123}t_{246})p_{3}^{2} \\ &+ (t_{124}t_{125} - t_{123}t_{126})(-p_{2}p_{3} + p_{1}p_{4}).
\end{aligned}
\end{equation*}}
\item We set 
\[
\Pi_{\mA}^{14}:=\Pi_{\mA}^{15}\cap\{t_{246}=1\}.
\]
\end{enumerate}
\end{defn}
\vspace{5pt}

For an alternative viewpoint on the equations $G_1,\dots,G_9$, the reader is referred to \cite[Rem.\,4.2]{Tak3}.

The equations $F_{1},\dots,F_{9}$ of $\Sigma_{\mA}^{13}$ with $s_{33}=1$,
and the equations of $G_{1},\dots,G_{9}$ of $\Pi_{\mA}^{14}$ with
$t_{246}=1$ are explicitly listed in the \textit{Mathematica} code
\cite{Tak8} as data.

\subsection{Properties of the key varieties}

The following results are obtained in \cite[Prop.\,2.3, Lem.\,2.4, Cor.\,2.6, Prop.\,2.12 and 2.\,13]{Tak9}
for $\Sigma_{\mA}^{14}$ and in \cite[Prop.\,5.5, Lem.\,5.6, Cor.\,5.7, Prop.\,5.10,  5.12, and 6.1]{Tak6}
for $\Pi_{\mA}^{15}$.
\begin{prop}
\label{prop:Sigmamain} The affine schemes $\Sigma_{\mA}^{14}$ and
$\Pi_{\mA}^{15}$ satisfy the following properties:

\begin{enumerate}[$(1)$]

\item The scheme $\Sigma_{\mA}^{14}$ and $\Pi_{\mA}^{15}$ are respectively
$14$- and $15$-dimensional irreducible normal varieties with only
factorial Gorenstein terminal singularities.

\item Let $\mA_{\overline{\Sigma}}$ (resp.~$\mA_{\overline{\Pi}}$)
be the affine space with the entries of $S$ and $\bm{t}$ (resp.~
\[
t_{1},t_{2},t_{123},t_{124},t_{125},t_{126},t_{135},t_{136},t_{245},t_{246}\text{{\it )}}
\]
as the coordinates. Let $\widehat{\Sigma}$ (resp.~$\widehat{\Pi}$)
be the variety obtained from $\Sigma_{\mA}^{14}$ (resp.~$\Pi_{\mA}^{15}$)
by considering the coordinates of $\mA_{\Sigma}$ (resp.~$\mA_{\Pi})$
except those of $\mA_{\overline{\Sigma}}$ (resp.~$\mA_{\overline{\Pi}}$)
as projective coordinates. The natural projection $\rho_{\Sigma}\colon\widehat{\Sigma}\to\mA_{\overline{\Sigma}}$
(resp.~$\rho_{\Pi}\colon\widehat{\Pi}\to\mA_{\overline{\Pi}}$) is
generically a $\mP^{2}\times\mP^{2}$-fibration with the relative
Picard number one.

\item Each of the defining ideals of $\Sigma_{\mA}^{14}$ and $\Pi_{\mA}^{15}$
is generated by $9$ elements with $16$ relations.

\end{enumerate}
\end{prop}
We also see that $\Sigma_{\mA}^{13}$ and $\Pi_{\mA}^{14}$ have similar
properties since the open subsets $\Sigma_{\mA}^{14}\cap\{s_{33}\not=0\}$
of $\Sigma_{\mA}^{14}$ and $\Pi_{\mA}^{15}\cap\{t_{246}\not=0\}$
of $\Pi_{\mA}^{15}$ are isomorphic to $\Sigma_{\mA}^{13}\times\mA^{1*}$
and $\Pi_{\mA}^{14}\times\mA^{1*}$ respectively. 

\begin{prop}
\label{prop:candiv}Let $\mathfrak{K}$ be $\Sigma_{\mP}^{12}$ or
$\Pi_{\mP}^{13}$, where $\Sigma_{\mP}^{12}$ and $\Pi_{\mP}^{13}$
are the weighted projectivization of $\text{\ensuremath{\Sigma_{\mA}^{13}}and }\Pi_{\mathbb{\mathbb{A}}}^{14}$
respectively with some positive weights of coordinates. For $\mathfrak{K}=\Sigma_{\mP}^{12}$,
we set 
\[
d_{0}:=w(p_{1})+w(q_{1})(=w(p_{2})+w(q_{2})=w(p_{3})+w(q_{3})).
\]
It holds that $\omega_{\mathfrak{K}}=\sO(-k)$, where 

$k=\begin{cases}
4w(r)+2w(u)+3(w(p_{1})+w(p_{2})+w(p_{3}))-6d_{0} & :\mathfrak{K}=\Sigma_{\mP}^{12}\\
7w(s_{3})-5w(p_{3})-2w(u_{1})+3w(u_{2}) & :\mathfrak{K}=\Pi_{\mP}^{13}.
\end{cases}$
\end{prop}
\begin{prop}
\label{prop:Pic1} Let $\mathfrak{K}$ be $\Sigma_{\mP}^{12}$ or
$\Pi_{\mP}^{13}$ as in Proposition \ref{prop:candiv}, and $\mathfrak{K}_{\mA}=\Sigma_{\mA}^{13}$
or $\Pi_{\mA}^{14}$ respectively. Let 
\[
\mathsf{b}:=\begin{cases}
p_{1} & :\mathfrak{K}=\Sigma_{\mP}^{12}\\
p_{3}^{2}+t_{1}p_{4}^{2}+p_{4}u_{2} & :\mathfrak{K}=\Pi_{\mP}^{13}.
\end{cases}
\]
It holds that

\begin{enumerate}[$(1)$]

\item $\mathfrak{K}_{\mA}\cap\{\mathsf{b}=0\}$ is irreducible and
normal, and the affine coordinate ring of $\mathfrak{K}_{\mA}\setminus\{\mathsf{b=0\}}$
is a UFD. This implies that any prime Weil divisor on $\mathfrak{K}$
is the restriction of a weighted hypersurface. In particular, $\mathfrak{K}$
is $\mathbb{Q}$-factorial and have Picard number one, and

\item let $X$ be a quasi-smooth $3$-fold such that $X$ is a codimension
$\dim\mathfrak{K}-3$ weighted complete intersection in $\mathfrak{K}$.
Assume moreover that $X\cap\{\mathsf{b}=0\}$ is a prime divisor.
Then any prime Weil divisor on $X$ is the restriction of a weighted
hypersurface. In particular, $X$ is $\mathbb{Q}$-factorial and has
Picard number one.

\end{enumerate}
\end{prop}
We remark that the proposition remains valid in the case $\mathfrak{K}=$
$\Pi_{\mP}^{14}$ (namely, $\mathfrak{K}_{\mA}=$the cone over $\Pi_{\mA}^{14})$,
and that the proof works verbatim.

\begin{acknowledgements} The author owes many important
calculations in the paper to Professor Shinobu Hosono. The author
wishes to thank him for his generous cooperations. The author would
like to express his sincere gratitude to the referee for their careful
reading of the manuscript and for offering many constructive suggestions.
This work is supported in part by Grant-in Aid for Scientific Research
(C) 16K05090. 
\end{acknowledgements}

\end{document}